\documentclass[journal]{IEEEtran}

\usepackage[T1]{fontenc}
\usepackage[utf8]{inputenc} 
\usepackage{graphicx}
\usepackage[]{hyperref}
\usepackage[version=3]{mhchem}
\usepackage{amssymb}
\usepackage{longtable}
\usepackage{multirow}
\usepackage{amsmath}
\usepackage{eurosym}	 
\usepackage{mathtools}
\usepackage{booktabs}

\usepackage[T1]{fontenc}
\usepackage[version=3]{mhchem}
\usepackage{amssymb}
\usepackage{longtable}
\usepackage{multirow}
\usepackage{amsmath}
\usepackage{mathtools}
\usepackage{epstopdf}
\usepackage{caption}
\usepackage{color}
\usepackage{xcolor}
\definecolor{myred}{rgb}{0.7,0.15,0.15}
\definecolor{mygreen}{rgb}{0.13,0.55,0.13}
\definecolor{myblue}{rgb}{0.25,0.41,0.88}

\usepackage{graphicx}
\usepackage[labelformat=simple]{subcaption} 

\usepackage{placeins}
\usepackage[font={small}]{caption}

\usepackage{mathrsfs}
\usepackage{comment}

\usepackage{enumitem}
\usepackage{url}
\usepackage{times}

\usepackage{pgfplots}
\usepackage{tikz}
\usetikzlibrary{patterns}
\usetikzlibrary{spy}

\usepackage{pgfkeys}
\usetikzlibrary{intersections}
\usepackage{listings}
\renewcommand{\vec}[1]{\boldsymbol{#1}} 
\newcommand{\mat}[1]{\mathbf{#1}} 

\newcommand{\paren}[1]{\left(#1\right)}
\newcommand{\parenangle}[1]{\left\langle#1\right\rangle}

\newcommand{\volInt}[3]{\paren{#1 \, , #2}_{#3}}
\newcommand{\volIntBig}[3]{\big(#1 \, , #2\big)_{#3}}

\newcommand{\surInt}[3]{\parenangle{#1 \, , #2}_{#3}}

\makeatletter
\newcommand\Label[1]{&\refstepcounter{equation}\text{(\theequation)}\ltx@label{#1}&}
\makeatother

\newlength\Colsep
\newcommand{\grad}{\text{\textbf{grad}}\, }

\renewcommand{\div}{\text{div}\, }

\newcommand{\curl}{\text{\textbf{curl}}\, }
\newcommand{\curlOnly}{\text{\textbf{curl}}}

\newcommand{\Def}{\coloneqq}

\renewcommand{\b}{\vec b}
\renewcommand{\a}{\vec a}
\renewcommand{\t}{\vec t}
\newcommand{\n}{\vec n}
\newcommand{\h}{\vec h}

\newcommand{\mur}{\mu_\text{r}}
\newcommand{\e}{\vec e}

\renewcommand{\j}{\vec j}

\newcommand{\dt}{\partial_t}
\newcommand{\ec}{e_{\text{c}}}
\newcommand{\jc}{j_{\text{c}}}

\renewcommand{\O}{\Omega}
\newcommand{\Oa}{\Omega_a}
\newcommand{\Oh}{\Omega_h}
\newcommand{\Ohc}{\Omega_{h,\text{c}}}
\newcommand{\Ohcc}{\Omega_{h,\text{c}}^{\text{C}}}
\newcommand{\Ohci}{\Omega_{h,\text{c}_i}}

\newcommand{\Gm}{\Gamma_\text{m}}
\newcommand{\af}{$a$-formulation\ }
\newcommand{\hf}{$h$-formulation\ }
\newcommand{\haf}{$h$-$a$-formulation\ }
\newcommand{\taf}{$t$-$a$-formulation\ }
\newcommand{\afOnly}{$a$-formulation}
\newcommand{\hfOnly}{$h$-formulation}
\newcommand{\hafOnly}{$h$-$a$-formulation}
\newcommand{\tafOnly}{$t$-$a$-formulation}

\newcommand{\Od}{\Omega^\delta}
\newcommand{\Oad}{\Omega_a^\delta}
\newcommand{\Ohd}{\Omega_h^\delta}
\newcommand{\Ohcd}{\Omega_{h,\text{c}}^\delta}

\newcommand{\Gd}{\Gamma^\delta}
\newcommand{\Gmd}{\Gamma_\text{m}^\delta}
\newcommand{\Gwd}{\Gamma_w^\delta}
\newcommand{\GGwmd}{\partial\Gamma_w^{-,\delta}}
\newcommand{\GGwpd}{\partial\Gamma_w^{+,\delta}}
\newcommand{\GGwpid}{\partial\Gamma_{w,i}^{+,\delta}}

\newcommand{\Odb}{\bar{\Omega}^\delta}
\newcommand{\Oadb}{\bar{\Omega}_a^\delta}

\newcommand{\Ohcdb}{\bar{\Omega}_{h,\text{c}}^\delta}
\newcommand{\Ohccdb}{\bar{\Omega}_{h,\text{c}}^{\text{C},\delta}}

\newcommand{\Gwdb}{\bar{\Gamma}_w^\delta}

\newcommand{\hsp}{\mathcal{H}}
\newcommand{\asp}{\mathcal{A}}
\newcommand{\tsp}{\mathcal{T}}

\newcommand{\hspz}{\mathcal{H}_{0}}
\newcommand{\aspz}{\mathcal{A}_{0}}
\newcommand{\tspz}{\mathcal{T}_{0}}

\newcommand{\hspzd}{\mathcal{H}_{0}^{\delta}}
\newcommand{\aspzd}{\mathcal{A}_{0}^{\delta}}
\newcommand{\tspzd}{\mathcal{T}_{0}^{\delta}}

\newcommand{\hspdone}{\mathcal{H}^{\delta,1}}
\newcommand{\aspdone}{\mathcal{A}^{\delta,1}}
\newcommand{\tspdone}{\mathcal{T}^{\delta,1}}
\newcommand{\hspdtwo}{\mathcal{H}^{\delta,2}}
\newcommand{\aspdtwo}{\mathcal{A}^{\delta,2}}
\newcommand{\tspdtwo}{\mathcal{T}^{\delta,2}}

\newcommand{\hspzdone}{\mathcal{H}_{0}^{\delta,1}}
\newcommand{\aspzdone}{\mathcal{A}_{0}^{\delta,1}}
\newcommand{\tspzdone}{\mathcal{T}_{0}^{\delta,1}}
\newcommand{\hspzdtwo}{\mathcal{H}_{0}^{\delta,2}}
\newcommand{\aspzdtwo}{\mathcal{A}_{0}^{\delta,2}}
\newcommand{\tspzdtwo}{\mathcal{T}_{0}^{\delta,2}}
\newcommand{\hspzdi}{\mathcal{H}_{0}^{\delta,i}}

\newcommand{\tspzdi}{\mathcal{T}_{0}^{\delta,i}}

\newcommand{\aspzdj}{\mathcal{A}_{0}^{\delta,j}}

\newcommand{\ad}{\vec a^{\delta}}
\newcommand{\bd}{\vec b^{\delta}}
\newcommand{\td}{\vec t^{\delta}}
\newcommand{\hd}{\vec h^{\delta}}
\newcommand{\jzd}{ j^{\delta}_z}

\newcommand{\transpose}{^{\text T}}

\usepackage{eso-pic}

\begin{document}

\title{On the Stability of Mixed Finite-Element Formulations for High-Temperature Superconductors}

\author{Julien~Dular, Mané~Harutyunyan, Lorenzo~Bortot, Sebastian~Schöps,
        Benoît~Vanderheyden,
        and~Christophe~Geuzaine
\thanks{J. Dular, B. Vanderheyden, and C. Geuzaine are with the Department of Electrical Engineering and Computer Science, Institut Montefiore B28 in University of Liège, B-4000 Liège, Belgium.}%
\thanks{M. Harutyunyan, L. Bortot, and S. Schöps are from Computational Electromagnetics, Technical University of Darmstadt, Germany.}
\thanks{L. Bortot is with CERN, Geneva, Switzerland.}%
\thanks{J. Dular is a research fellow funded by the F.R.S-FNRS.}}

\maketitle

\AddToShipoutPicture*{
    \footnotesize\sffamily\raisebox{0.8cm}{\hspace{1.5cm}\fbox{
        \parbox{\textwidth}{
            \copyright~2021
                IEEE. Personal use of this material is permitted. Permission from IEEE
                must be obtained for all other uses, in any current or future media, 
                including reprinting/republishing this material for advertising or
                promotional purposes, creating new collective works, for resale or
                redistribution to servers or lists, or reuse of any copyrighted
                component of this work in other works.
            }
        }
    }
}

\begin{abstract}

In this work, we present and analyze the numerical stability of two coupled finite element formulations. The first one is the \haf and is well suited for modeling systems with superconductors and ferromagnetic materials. The second one, the so-called \taf with thin-shell approximation, applies for systems with thin superconducting domains. Both formulations involve two coupled unknown fields and are mixed on the coupling interfaces. Function spaces in mixed formulations must satisfy compatibility conditions to ensure stability of the problem and reliability of the numerical solution. We propose stable choices of function spaces using hierarchical basis functions and demonstrate the effectiveness of the approach on simple 2D examples.
\end{abstract}


\IEEEpeerreviewmaketitle


\section{Introduction}
\addcontentsline{toc}{section}{Introduction}

\IEEEPARstart{M}{odeling} accurately and efficiently the magnetic response of high-temperature superconductors (HTS) is important for the development of numerous magnet and electrical power applications, e.g., superconducting rotating machines. One of the main tools used to model the properties of superconductors is the finite element method (FEM), based on formulations of Maxwell's equations combined with the $E$-$J$ power law. This law is strongly non-linear and requires a carefully chosen formulation.


In the past few years, several FEM models based on coupled formulations have been proposed. An \haf was introduced in a 2D model of rotating machines with superconducting windings~\cite{Brambilla2018}. The superconducting materials were modeled with an \hfOnly, whereas the continuity conditions between the fixed and the rotating parts were treated with the \afOnly. A second \haf formulation was introduced for modeling HTS magnets with a coupling to an external circuit, in order to reduce the number of degrees of freedom in the non-superconducting regions~\cite{Bortot} with respect to a full \hf~\cite{Hong2006,Shen2020}. A third \haf was considered for systems containing superconductors and ferromagnets, in order to model each material with its most efficient formulation~\cite{dular2019finite}. A setting involving the simultaneous computation of magnetic and electric fields in the whole conducting domain has been proposed in \cite{barrett2012electric} for thin superconducting films. Another type of combined formulation, the \tafOnly, was introduced in \cite{zhang2016efficient} to model superconducting tapes, presenting a high width over thickness ratio. The current density inside the tapes was described by a surface current potential, whereas the magnetic vector potential was the state variable outside the tape. Both fields were coupled by means of integrals on the surface of the (infinitely thin) HTS tape. In~\cite{Bortot}, the same formulation is derived from the \haf with a thin-sheet approximation. The \taf has also recently been extended to finite volume systems, e.g., by modeling a stack of tapes in full or in parts as an equivalent homogeneous bulk material~\cite{berrospe2019real,Wang2020}.

In each of these coupled formulations, different finite element fields
are introduced region-wise, while they coexist and are coupled through
a common boundary or a common region. The coupling makes these
formulations \textit{mixed}, for which care must be taken in the
choice of function spaces and the discretization. For instance, naive
choices of approximation function spaces can easily lead to stability
issues manifesting themselves as spurious oscillations in the numerical
solution (chapter 8 of Ref.~\cite{brezziBook}). Such oscillations have been indeed
observed numerically in the \taf \cite{berrospe2019real} when
using first-order polynomials for both the $t$- and $a$-approximation
spaces. General mathematical conditions for solvability and stability
have been stated and studied in a number of mixed finite element
problems ~\cite{brezziBook, bathe2001inf, babuvska1973finite}, both for the continuous
and the discrete problems. One of these conditions, known as the
inf-sup condition, is usually difficult to prove analytically but may
be tested numerically~\cite{bathe2001inf}.

The problems we consider are nonlinear but the stability issues and the resulting
oscillations are not a direct consequence of the nonlinearity of the
constitutive laws. We observed that they actually already appear in linear problems
with the same coupled formulations and non-compatible function
spaces. 
However, the nonlinearity of the constitutive
laws is one of the motivations for using coupled
formulations, which is the reason why we discuss their stability in the framework of superconducting systems.

In this work, we consider the \haf for systems
containing superconductors and ferromagnets and the \tafOnly, which can be seen as the limit of the \haf for thin superconducting tapes. Following the general theory of mixed finite
elements, we analyze the related conditions for obtaining numerically
stable mixed formulations. In section \ref{sec_formulations}, we
introduce and derive the two coupled formulations. In particular, we derive a
version of the \taf that directly includes global variables on current
intensity or voltage in the weak form. To the best of our knowledge, it
has not been introduced in that form yet. In section
\ref{sec_discretization}, we present several choices of discretized
function spaces and illustrate the occurrence of the spurious
oscillations that arise when spaces are not compatible. We
recall the classical theory of mixed formulations and perturbed
saddle-point problems \cite{brezziBook} in section
\ref{sec_stabilityAnalysis} and present a \textit{numerical inf-sup
test} based on \cite{bathe2001inf} to check the compatibility of discretized function spaces.
In the last two sections, the theory is applied to the \haf (section
\ref{sec_haf}) and the \taf (section \ref{sec_taf}), restricting to 2D problems with in-plane magnetic fields.

\section{Finite-Element Formulations}\label{sec_formulations}

The magnetic response of a system containing type-II superconductors with strong pinning can be described by Maxwell's equations in the magnetodynamic (quasistatic) approximation \cite{jackson1999classical}, and magnetic and electrical constitutive laws,
\begin{equation}\label{MQSequations}
\left\{\begin{aligned}
\div\b &= 0,& \text{(magn. Gauss)}\\
\curl\h &= \j, &\text{(Ampère)}\\
\curl\e &= -\dt \b,&\text{(Faraday)}
\end{aligned}\right. \quad \text{and} \quad  \left\{\begin{aligned}
\b &= \mu \h,\\
\e &= \rho \j,
\end{aligned}\right.
\end{equation}
with $\b$, $\h$, $\j$, $\e$, $\mu$, and $\rho$, being the magnetic flux density (T), the magnetic field (A/m), the electric current density (A/m$^2$), the electric field (V/m), the permeability (H/m), and the resistivity ($\O$m), respectively. The permeability can be a function of $\h$. In non-conducting materials, $\rho \to \infty$ and $\j = \vec 0$. In superconductors, $\b = \mu_0 \h$ and we assume a power law for the resistivity \cite{rhyner1993magnetic},
\begin{equation}\label{eqn_contitutiveje}
\rho = \frac{\ec}{\jc}\paren{ \frac{\|\j\|}{\jc}}^{n-1},
\end{equation}
where $\ec = 10^{-4}$ V/m is a threshold electric field defining the critical current density $\jc$ (A/m$^2$). The dimensionless number $n=U_0/k_\text{B}T$, with $U_0$ a pinning energy and $T$ the temperature, is a critical exponent associated with magnetic flux creep.

In the following, the system is modeled in a domain $\O$. Boundary conditions are applied on its external boundary $\partial \O$, which is decomposed into two complementary domains: $\Gamma_e$, where the normal component of $\b$ or the tangential component of $\e$ is imposed, and $\Gamma_h$, where the tangential component of $\h$ is imposed. We also use the following notation for volume and surface integrals:
\begin{align}
\volInt{f_1}{f_2}{\Omega} = \int_{\Omega} f_1 \cdot f_2 \ d\Omega, \quad  \surInt{f_1}{f_2}{\Gamma} = \int_{\Gamma} f_1 \cdot f_2 \ d\Gamma,
\end{align}
with $f_1$ and $f_2$ being two scalar or vector fields and $\cdot$ the scalar multiplication or the dot product, respectively.

We now present two mixed finite element formulations of the magnetodynamic problem. 

\subsection{Coupled formulation 1 - \haf}

When a system contains both a superconductor and a nonlinear ferromagnetic material, classical formulations such as the \hf or the \af may face convergence issues. The power law in superconductors is easier to handle with a Newton-Raphson method in the \hfOnly, which involves the electrical resistivity. 
Conversely, the \af is more efficient than the \hf to deal with the typical saturation law describing the permeability of ferromagnets~\cite{dular2019finite}. Combining the \hf and \af into a coupled \haf by choosing the best formulation in each region has proven to be an efficient solution for systems with both materials~\cite{dular2019finite}.

\begin{figure}[h!]
\centering
            \begin{subfigure}[b]{0.49\linewidth}
            \centering
		\includegraphics[width=\textwidth]{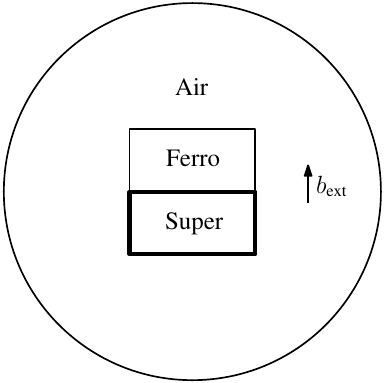}
		\caption{Problem geometry.}
		\label{bar_mesh}
        \end{subfigure}
\begin{subfigure}[b]{0.49\linewidth}  
            \centering 
		\includegraphics[width=\textwidth]{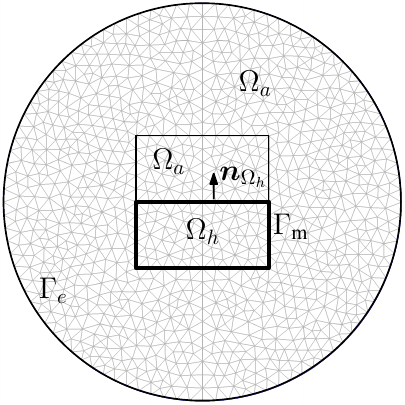}
		\caption{Domain definition and mesh.}
		\label{bar_a1_h1}	
      \end{subfigure}
        \caption{Stack of a superconducting bar (below,  $n=20$, $\jc=3\times 10^8$ A/m$^2$) and a ferromagnetic bar (above, $\mu_\text{r} = 1000$), subjected to an external field ($b_{\text{ext}} = 0.4$ T). The thick curve is $\Gm$.}
        \label{bar_b}
\end{figure}

The domain $\O$ is decomposed into two parts: $\Oh$, containing the superconducting domain, and $\Oa$, containing the nonlinear ferromagnetic domain, which is assumed to have a negligible electrical conductivity. The parts of $\O$ where constitutive laws are linear can be put in $\Oh$ or $\Oa$. Inside $\Oh$, the conducting domain is denoted by $\Ohc$, and the non-conducting domain is denoted by $\Ohcc$, with $\Oh = \Ohc \cup \Ohcc$. The common boundary of $\Oh$ and $\Oa$ is denoted by $\Gm$. Coupling operates via this common interface. We also introduce the outer normal vectors $\n_{\Oh}$ and $\n_{\Oa}$. For illustration, consider the simple 2D stacked bar geometry in Fig.~\ref{bar_b}, where the \hf is applied to a superconducting bar, the \af is applied to a ferromagnetic bar and to the air region, while the coupling surface $\Gm$ is taken as the boundary of the superconducting region. In this example, the external boundary belongs to $\Gamma_e$, and $\Gm$ is the boundary of the superconducting region, which constitutes the entire $\Oh$ domain.

We derive the two formulations in $\Oa$ and $\Oh$, then couple them to obtain the \hafOnly.

The classical \af is a weak form of Ampère's law where the magnetic flux density $\b$ is expressed via a vector potential $\a$ as $\b = \curl \a$. Here, it is introduced in $\Oa$ only. We choose $\a\in \asp(\Oa)$ with 
\begin{align}\label{eqn_asp}
\asp(\Oa) &= \big\{\a \in H(\curlOnly;\Oa)\ |\notag \\
&\qquad \quad \a \times \n_{\Oa} =  \bar \a \times \n_{\Oa} \text{ on } \Gamma_e\cap \partial \Oa\big\},
\end{align}
with $\bar \a\times \n_{\Oa}$ a fixed trace on $\Gamma_e\cap \partial \Oa$. For conciseness, we place all conducting materials in $\Oh$, and Ampère's law reads $\curl \h = \vec 0$ in $\Oa$. We multiply this equation by a test function $\a'$ in the space $\aspz(\Oa)$ with homogeneous essential boundary conditions $\a\times \n_{\Oa}=\vec 0$ on $\Gamma_e \cap \partial \Oa$,
and integrate the product over $\Oa$. We obtain
\begin{align}
& \volInt{\curl \h}{\a'}{\Oa} = 0 \notag\\
\Leftrightarrow \ & \volInt{\h}{\curl \a'}{\Oa} - \surInt{\h\times \vec n_{\Oa}}{\a'}{(\Gamma_h\cap \partial \Oa) \cup \Gm} = 0,
\end{align}
using a Green identity. Prescribing the value of $\h\times \n_{\Oa}$ on $\Gamma_h\cap \partial \Oa$ constitutes a natural boundary condition for the \afOnly. For conciseness, we consider homogeneous natural boundary conditions on $\Gamma_h\cap \partial \Oa$. Therefore, after introducing the vector potential $\a$, the formulation amounts to finding $\a\in \asp(\Oa)$ such that $\forall \a'\in \aspz(\Oa)$,
\begin{align}\label{eqn_aformulation_ha}
\volInt{\nu\,\curl \a}{\curl \a'}{\Oa} - \surInt{\h\times \vec n_{\Oa}}{\a'}{\Gm} = 0,
\end{align}
with the reluctivity $\nu = \mu^{-1}$. On $\Gm$, the tangential magnetic field is still unknown. It will be coupled with the formulation in $\Oh$ that we derive next.

In $\Oh$, we use the \hf with curl-free functions in $\Ohcc$, also called the $h$-$\phi$-formulation. This is a weak form of Faraday's law. We consider $N$ distinct conducting subdomains $\Ohci$ of $\Ohc$, with $i\in C=\{1,2,\dots,N\}$, on which we impose either the current or the voltage. The current is imposed on a subset $C_I$ of $C$, and the voltage is imposed on the complementary set $C_V$. In the \hfOnly, the magnetic field $\h$ is sought in $\hsp(\Oh)$ defined as
\begin{align}\label{eqn_hsp}
&\hsp(\Oh) = \big\{\h \in H(\curlOnly;\Oh)\ |\ \curl \h = \vec 0 \text{ in } \Ohcc,\notag \\
&\quad \h \times \n = \bar \h \times \n_{\Oh} \text{ on } \Gamma_h\cap \partial \Oh,\, \mathcal{I}_i(\h) = I_i \text{ for } i\in C_I\big\}.
\end{align}
Only curl-free functions are considered for $\h$ in the non-conducting domain $\Ohcc$, so that the current density $\j = \curl \h$ is exactly zero in $\Ohcc$, by construction. Functions associated with net electrical currents in the conducting regions (that are not gradients of a scalar potential) are however still considered, they form a cohomology basis of dimension $N$ \cite{pellikka2013homology}. Each of these functions can be associated with a (group of) conducting subdomain(s) $\Ohci$ of $\Ohc$, $i\in C=\{1,2,\dots,N\}$. The notation $\mathcal{I}_i(\h)$ denotes the net current $I_i$ flowing in (a group) of conductor(s) $i$ for a given function $\h$ \cite{dular1994phd}, i.e., the circulation of $\h$ along a closed loop $\mathcal{C}_i$ around that (group of) conductor(s):
\begin{align}
\mathcal{I}_i(\h) = \oint_{\mathcal{C}_i} \h \cdot d\vec \ell.
\end{align}
The corresponding applied voltage is denoted below by $V_i$ (voltage per unit length in 2D). Either $I_i$ or $V_i$ must be imposed for each $i$. Note that imposing a zero curl in $\Ohcc$ drastically reduces the number of degrees of freedom in the function space \eqref{eqn_hsp}, after discretization. Indeed, in 2D problems with in-plane magnetic field or 3D problems, with Whitney basis functions \cite{bossavit1988whitney}, only one unknown per node is necessary, compared to one unknown per edge for the "full \hfOnly". 

The weak form is obtained by projecting Faraday's law on test functions, $\forall \h'\in \hspz(\Oh)$,
\begin{align}
& \volInt{\dt(\mu\, \h)}{\h'}{\Oh} + \volInt{\curl \e}{\h'}{\Oh} = 0\notag\\
\Leftrightarrow\ & \volInt{\dt(\mu\, \h)}{\h'}{\Oh} + \volInt{\e}{\curl \h'}{\Oh}\notag \\
&\qquad\qquad-\surInt{\e \times \n_{\Oh}}{\h'}{(\Gamma_e\cap \partial \Oh)\cup \Gm}= 0.
\label{eqn_hformulation_ha_derivation}
\end{align}
The space $\hspz(\Oh)$ for test functions is with homogeneous essential boundary conditions, $\h \times \n_{\Oh} = \vec 0$ on $\Gamma_h\cap \partial \Oh$ and $\mathcal{I}_i(\h) = 0$ for $i\in C_I$.

In $\Ohcc$, $\curl \h'=\vec 0$. In $\Ohc$, $\e = \rho\, \curl \h$. We model localized power sources and their associated current and voltage on each conducting subdomain as proposed in \cite{dular1994phd,dular1999global}. For conciseness again, we consider homogeneous natural boundary conditions on $\Gamma_e\cap \partial \Oh$. Formulation \eqref{eqn_hformulation_ha_derivation} then becomes
\begin{align}\label{eqn_hformulation_ha}
& \volInt{\dt(\mu\, \h)}{\h'}{\Oh} + \volInt{\rho\,\curl \h}{\curl \h'}{\Ohc}\notag \\
&\qquad\qquad-\surInt{\e \times \n_{\Oh}}{\h'}{\Gm}=  -\sum_{i\in C} V_i \mathcal{I}_i(\h'),
\end{align}
with the $V_i$'s being natural "boundary" conditions for $i\in C_V$. For $i\in C_I$, the current $I_i$ is imposed, then $\mathcal{I}_i(\h') = 0$ and the global term does not enter the problem. It can however be exploited to build a circuit equation, to compute the voltage $V_i$ associated with the imposed current $I_i$ as a post-processing quantity. Conversely, for $i\in C_V$, the voltage $V_i$ is imposed, then $I_i$ is a degree of freedom and the global term enters the system of equations \cite{dular1999global}.

On $\Gm$, the tangential electric field $\e$ is still unknown. It will be coupled with the formulation in $\Oa$.

The final step in the \haf derivation consists in coupling the two separate formulations \eqref{eqn_aformulation_ha} and \eqref{eqn_hformulation_ha} in $\Oa$ and $\Oh$. The tangential trace of the magnetic field on $\Gm$ in \eqref{eqn_aformulation_ha} can be directly expressed in terms of the magnetic field $\h$ of \eqref{eqn_hformulation_ha}. Similarly, the tangential trace of the electric field on $\Gm$ in \eqref{eqn_hformulation_ha} can be expressed in terms of the vector potential $\a$ of \eqref{eqn_aformulation_ha}, with $\e = -\dt \a - \grad v$. In fact, only the $-\dt \a$ term contributes to the integral (see Appendix).

The resulting coupled \haf reads:

From an initial solution at time $t=0$, find $\h \in \hsp(\Oh)$ and $\a\in \asp(\Oa)$ such that, for $t>0$, $\forall \h' \in \hspz(\Oh)$, and $\forall \a' \in \aspz(\Oa)$,
\begin{equation}\label{eqn_coupledFormulation}
\begin{aligned}
& \volIntBig{\dt (\mu\, \h)}{\h'}{\Oh} + \volIntBig{\rho\,  \curl \h}{\curl \h'}{\Ohc}\\
& \qquad  + \surInt{\dt \a \times\vec n_{\Oh}}{\h'}{\Gamma_\text{m}}= -\sum_{i\in C} V_i \mathcal{I}_i(\h'),\\
&\surInt{\h\times\vec n_{\Oa}}{\a'}{\Gamma_\text{m}} - \volInt{\nu\, \curl \a}{\curl \a'}{\Oa} = 0.
\end{aligned}
\end{equation}

The discrete function spaces must be chosen with care. In particular, the choice of basis functions spanning the trace space on $\Gm$ will affect the stability of the method. Different possibilities will be analyzed in section~\ref{sec_discretization}.

\subsection{Coupled formulation 2 - \taf for thin tapes}

The second formulation we consider is the so-called \taf for modeling thin superconducting tapes \cite{zhang2016efficient}. The tape is modeled as a line in 2D (a surface in 3D). The current density inside the tape is described via a current vector potential whereas the external magnetic flux density is expressed as the curl of a magnetic vector potential, naturally allowing discontinuous tangential components of the magnetic field across the tape. This \taf can be viewed as a geometric  limit case of the \haf for thin geometries. It possesses the same continuity properties and advantages.

In \cite{Bortot}, this formulation is derived from the \haf with a thin-sheet approximation. Circuit coupling is then performed by means of winding functions \cite{schops2013winding}.

Here, we present a version of the \taf following a different approach for circuit coupling. With the same philosophy as in \cite{dular1999global}, in each tape we either strongly impose the current intensity, directly in the function space, or weakly impose the voltage, with a circuit equation contained in the formulation. The formulation is valid in 2D or 3D. The stability analysis will be conducted in 2D in section~\ref{sec_stabilityAnalysis}.

\begin{figure}[h!]
            \begin{subfigure}[b]{0.49\linewidth}
            \centering
		\includegraphics[width=\textwidth]{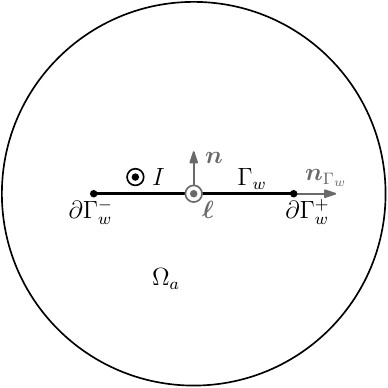}
		\caption{2D problem.}
		\label{tapeGeometry_ta_article_2D}
        \end{subfigure}
\begin{subfigure}[b]{0.49\linewidth}  
            \centering 
		\includegraphics[width=\textwidth]{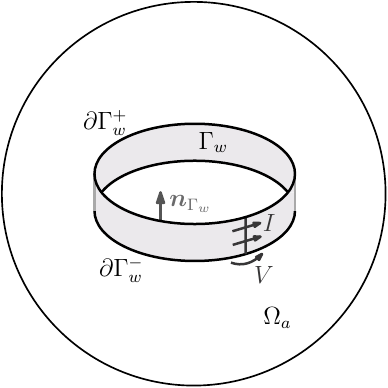}
		\caption{3D problem (for illustration).}
		\label{tapeGeometry_ta_article_3D}	
      \end{subfigure}
\caption{Conventions for the \taf derivation. (a) 2D case, a tape with current density perpendicular to the modeled plane. (b) 3D case with a tape loop, e.g., a racetrack coil. In 3D, the effect of an external voltage/current source is modeled on an arbitrary cross-section.}
\label{tapeGeometry_ta}
\end{figure}

This \taf applies to situations with thin conducting domains. Let us consider a conducting domain $\Gamma_w\subset \Oa$ of thickness $w$, see Fig.~\ref{tapeGeometry_ta}. We start from the classical \af in the whole domain $\Oa$, with homogeneous natural boundary conditions for conciseness: find $\a\in \asp(\Oa)$ such that $\forall \a' \in \aspz(\Oa)$,
\begin{equation}\label{a_formulation_ta_js}
\volInt{\nu\, \curl \a}{\curl \a'}{\Oa} - \volInt{\j}{\a'}{\Gamma_w} = 0,
\end{equation}
with a given current density $\j$ (A/m$^2$) in $\Gamma_w$. Instead of modeling the tape $\Gamma_w$ as a volume, we collapse it into a surface and replace $\j$ by a surface current density $\vec k = w\j$ (A/m), perpendicular to the normal vector $\n$. This constitutes the main approximation of the formulation: the thickness is not represented in the geometry but introduced inside the equation. Possible variations of $\j$ across the thickness are therefore chosen not to be modeled.

Definition~\eqref{eqn_asp} implies that the vector potential $\a\in \asp(\Oa)$ is continuous across $\Gamma_w$, but allows $\h \times \n = \nu \,\curl\a \times \n$ to be discontinuous. Actually, we can show that $(\h_1-\h_2)\times \vec n = \vec k$ is weakly satisfied (with an upward normal, $\h_1$ is the field on the top of the tape and $\h_2$ is the field below).


If the current density were known, the problem would be closed. Here, we want to represent eddy currents and an equation for the distribution of $\vec k$ is required. Since the current density is divergence free (magnetodynamic regime), we can express the current density $\j$ via a current vector potential $\vec t$ defined up to a gradient such that $\j = \curl \vec t$. To gauge $\vec t$, we choose it along the normal to the tape, i.e., $\t = t\n$ \cite{zhang2016efficient}.

For simplicity, in 3D, we restrict ourselves to closed current loops. The tape boundary $\partial \Gamma_w$ is decomposed into two disjoint parts, $\partial\Gamma_w^-$ and $\partial\Gamma_w^+$, as represented in Fig.~\ref{tapeGeometry_ta}. We model a possible power source on an arbitrary cross-section of the tape that imposes either a current intensity $I$ or a voltage $V$.
On lateral edges $\partial\Gamma_w^-$ and $\partial\Gamma_w^+$, $\j \cdot \n_{\Gamma_w} = 0$ so $\t$ is constant. Let us (strongly) fix it to $0$ on $\partial\Gamma_w^-$ and let its value, denoted by $T$, remain free on the other lateral edge $\partial\Gamma_w^+$. The value of $T$ is related to the total injected current intensity $I$. Indeed, on any cross-section $S$ of the tape, using Stokes' theorem,
\begin{align}
I &= \int_S \vec j \cdot d\vec S = \int_S \curl \vec t \cdot d\vec S = \oint_{\partial S} \vec t \cdot d\vec \ell_{\partial S} \notag \\
&= w(t|_{\partial\Gamma_w^+} - t|_{\partial\Gamma_w^-}) = wT.
\end{align}

To obtain a weak formulation for $\vec t$, we use Faraday's law, $\dt \b + \curl \e = \vec 0$, more specifically its component along $\n$. It amounts to finding $\t\in \tsp(\Gamma_w)$, such that $\forall \t'\in \tspz(\Gamma_w)$,
\begin{align}
& 0 = \surInt{\dt \b}{\t'}{\Gamma_w} + \surInt{\curl(\rho\, \j)}{\t'}{\Gamma_w}\notag \\
\Leftrightarrow\ & 0 = \surInt{\dt (\curl \a)}{\t'}{\Gamma_w} + \surInt{\curl(\rho\, \curl \t)}{\t'}{\Gamma_w}\notag \\
\Leftrightarrow\ & 0 = \surInt{\dt \a}{\curl \t'}{\Gamma_w} + \surInt{\rho\, \curl \t}{\curl \t'}{\Gamma_w}\notag \\
&\qquad - \surInt{(\e +\dt \a) \times \n_{\Gamma_w}}{\t'}{\partial\Gamma_w},\label{eqnToSimplify}
\end{align}
where we expressed the normal flux density $\b\cdot\n$ via the vector potential of the \afOnly.
Note that the outer normal $\n_{\Gamma_w}$ of $\Gamma_w$ arising from Green's identities is different from $\n$, see Fig.~\ref{tapeGeometry_ta}. 
The spaces $\tsp$ and $\tspz$ will be defined later.

The last term in \eqref{eqnToSimplify} is exploited to impose global quantities, such as current intensity or voltage.  The electric field in an \af is expressed as $\e = -\dt \a - \grad v$, with a scalar electric potential $v$. 
Because $\t' = \vec 0$ on $\partial \Gamma_n^-$, the line integral in \eqref{eqnToSimplify} reads
\begin{align}\label{globalTAterm}
&\surInt{(\e+\dt \a) \times \n_{\Gamma_w}}{\t'}{\partial\Gamma_w} = - \surInt{\grad v \times \n_{\Gamma_w}}{T'\n}{\partial\Gamma_w^+}\notag \\
&\qquad = - \surInt{\grad v \cdot \vec \ell}{T'}{\partial\Gamma_w^+} = -VT',
\end{align}
with $\vec \ell = \n_{\Gamma_w} \times \n$, and $V$ being the net potential difference (V) applied by the generator in 3D.  In 2D, $V$ is a voltage per unit length (V/m) in the out-of-plane direction. In the tape, either the total current $I$ or the associated voltage $V$ must be imposed. As with the \hafOnly, if the current $I$ is imposed, then $\t' = \vec 0$ on $\partial \Gamma_n^+$, and the equation does not enter the problem. It can however be used as a circuit equation to compute the voltage $V$ associated with the imposed current $I$, as a post-processing quantity. Conversely, if the voltage $V$ is imposed, then $I$ is a degree of freedom and Eq.~\ref{globalTAterm} enters the system of equations.

We now consider $N$ distinct tapes $\Gamma_{w,i}$ with $i\in C = \{1,2,\dots,N\}$. The union of these tapes is $\Gamma_w$. Current is imposed on a subset $C_I$ of $C$ whereas voltage is imposed on the complementary set $C_V$. For conciseness, we consider homogeneous natural boundary conditions on $\Gamma_h$. 
The \taf reads as follows:

From an initial solution, find $\a \in \asp(\Oa)$ and $\t \in \tsp(\Gamma_w)$, such that for all time instants and $\forall \a' \in \aspz(\Oa)$, $\forall \t' \in \tspz(\Gamma_w)$,
\begin{align}\label{taformulation}
& \volInt{\nu\, \curl \a}{\curl \a'}{\Oa} - \surInt{w\, \curl \t}{\a'}{\Gamma_w} = 0,\notag\\
&\begin{aligned}[t]
\surInt{w\, \dt \a}{\curl \t'}{\Gamma_w} + \surInt{w\, \rho\, \curl \t}{\curl \t'}{\Gamma_w} \\
 =-  \sum_{i\in C} V_i\mathcal{I}_i(\t'),
\end{aligned}
\end{align}
with $\mathcal{I}_i(\t') = w T_i'=I_i'$ being the net current flowing in tape $i$ for the potential $\t'$. The space $\tsp(\Gamma_w)$ (resp. $\tspz(\Gamma_w)$) is the set of functions $\t = t\n$ such that $\curl \t$ is in the dual space of the relevant trace space on $\Gamma_w$ of functions in $\asp(\Oa)$, with $\t = \vec 0$ on $\partial\Gamma_w^-$, and $\t = (I_i/w) \n$ (resp. $\t =\vec 0$) on $\partial\Gamma_{w,i}^+$ for $i\in C_I$. Since in 2D the vector potential $\a$ has only one out-of-plane component, $\asp(\Oa)$ can be identified with $H^1(\Oa)$ (see chapter 2 of Ref.~\cite{brezziBook}). With $\vec z$ being the direction of the current density, perpendicular to the 2D plane, if $\Gamma_w\cap \Gamma_e=\emptyset$ \cite{bechet2009stable}, we can choose $\t$ in
\begin{align}
&\tsp(\Gamma_w) = \big\{\t = t\n\ \big|\ (\vec z\cdot \curl \t) \in H^{-1/2}(\Gamma_w),\notag \\
&\quad t=0 \text{ on } \partial \Gamma_w^-\ , t =I_i/w \text { on } \partial \Gamma_{w,i}^+ \text{ for } i\in C_I \big\}.
\end{align}

As with the coupled $h$-$a$-formulation, the discrete function spaces must be chosen with care. In particular, the choice of basis functions spanning the trace space on $\Gamma_w$ will affect the stability of the method. Different possibilities will be analyzed in section~\ref{sec_discretization}.

\section{Discretization and Oscillations}\label{sec_discretization}

To proceed, we discuss different discretization schemes and their consequences on the stability of the coupled formulations.

For the numerical resolution, the domain $\O$ is discretized as $\O^\delta$ with a finite element mesh of characteristic size $\delta$. Function spaces for $\h$, $\a$ and $\t$ are approximated by basis functions on the finite elements and we denote the approximated functions by $\hd$, $\ad$, and $\td$. We then integrate over time with an implicit Euler method and solutions to nonlinear systems are obtained by Newton-Raphson iterations.

We focus on 2D problems, such as those represented in Figs.~\ref{bar_mesh} and \ref{tapeGeometry_ta_article_2D}. Finite element modeling is performed by GetDP \cite{getdp} and finite element meshes are generated by Gmsh \cite{gmsh}\footnote{Model files for the main test cases are available on \url{www.life-hts.uliege.be}.}.

All three fields $\hd$, $\ad$, and $\td$ of the coupled formulations \eqref{eqn_coupledFormulation} and \eqref{taformulation} are approximated by 1-forms \cite{lindell2004differential}. The \haf is $b$-conform in $\Oad$ and $h$-conform in $\Ohd$. The \taf is $b$-conform in $\Oad$ and the current density $\j=\curl \t$ in $\Gwd$ is a 2-form so that the continuity of its normal component is satisfied \cite{lindell2004differential}. Note that the lack of $h$-conformity for the \taf in $\Oad$ naturally allows the tangential magnetic field $\n\times( \nu\ \curl \ad \times \n)$ to be discontinuous across each tape, while the discontinuity strength is enforced weakly by means of the surface terms.

\subsection{Lowest order Whitney basis functions}

The simplest approximation spaces are generated by lowest order Whitney edge functions for the three fields \cite{bossavit1988whitney}. We use the following notation: $n \in \Odb$ or $e \in \Odb$ refers to nodes $n$ or edges $e$ in $\Od$ and on its boundary $\partial \Od$. To exclude entities on a boundary $\Gd$, we note $n$ (or $e$) $\in \Odb \setminus \Gd$ explicitly.

We build the magnetic field $\hd$ in $\Ohd$ as follows,
\begin{equation}\label{eqn_h_decomposition}
\hd = \sum_{e \in \Ohcdb\setminus \partial \Ohcd} h_e \, \vec \psi_e + \sum_{n \in \Ohccdb} \phi_n \, \grad \psi_n + \sum_{i \in C} I_i \, \vec c_i,
\end{equation}
with $\vec \psi_e$ being the edge function of edge $e$, $\psi_n$ the node function of node $n$, and $\vec c_i$ a discontinuous basis function associated with the cut related to conducting region $i$, defined on a transition layer. Note that $\grad \psi_n$ and $\vec c_i$ can be expressed as sums of edge functions \cite{dular1999global, dular1994phd}. We denote by $\hspdone(\Ohd)$ the space generated by these lowest order functions, including essential boundary conditions. We have $\hspdone(\Ohd) \subset \hsp(\Ohd)$. Coefficients $h_e$, $\phi_n$ and $I_i$ are the degrees of freedom for $\hd$.  We have $\mathcal{I}_i(\h^\delta) = I_i$ with the notation of formulation \eqref{eqn_coupledFormulation}, i.e., $I_i$ is the net current intensity flowing in (a group of) conductor(s) $i$ for the field $\h^\delta$.

The magnetic vector potential $\ad$ in $\Oad$ in both $h$-$a$ and $t$-$a$ formulations reads 
\begin{align}\label{eqn_a_decomposition}
\ad = \sum_{n \in \Oadb} a_n \, \psi_n\vec z,
\end{align}
where $\psi_n \vec z$ is a "perpendicular edge function" associated with node $n$, such that $\ad$ is chosen along  $\vec z$, the direction perpendicular to the plane in 2D. The field $\ad$ automatically satisfies the Coulomb gauge condition $\div\ad = 0$. 
We denote by $\aspdone(\Oad)$ the space generated by these lowest order functions, including essential boundary conditions. We also have $\aspdone(\Oad) \subset \asp(\Oad)$. Coefficients $a_n$ are the degrees of freedom for $\ad$. 

The current vector potential $\td$ reads
\begin{equation}\label{eqn_t_decomposition}
\td = \sum_{n\in \Gwdb\setminus (\GGwmd \cup \GGwpd)} t_n\, \psi_n\n + \sum_{i\in C} T_i\, \vec \ell_i,
\end{equation}
with $\psi_n \n$ being a "perpendicular edge function" associated with node $n$ and $\vec \ell_i$ the perpendicular edge function associated with the node on $\GGwpid$. (Note that in 3D, $\vec \ell_i$ is the sum of all perpendicular edge functions associated with nodes on $\GGwpid$; this is a global basis function.) We denote by $\tspdone(\Gwd)$ the space generated by these functions, including essential boundary conditions. Again, we have $\tspdone(\Gwd)\subset \tsp(\Gwd)$. Coefficients $t_n$ and $T_i$ are the degrees of freedom for $\td$.

\begin{figure}[h!]
\centering
            \begin{subfigure}[b]{\linewidth}
            \centering
		\includegraphics[width=\textwidth]{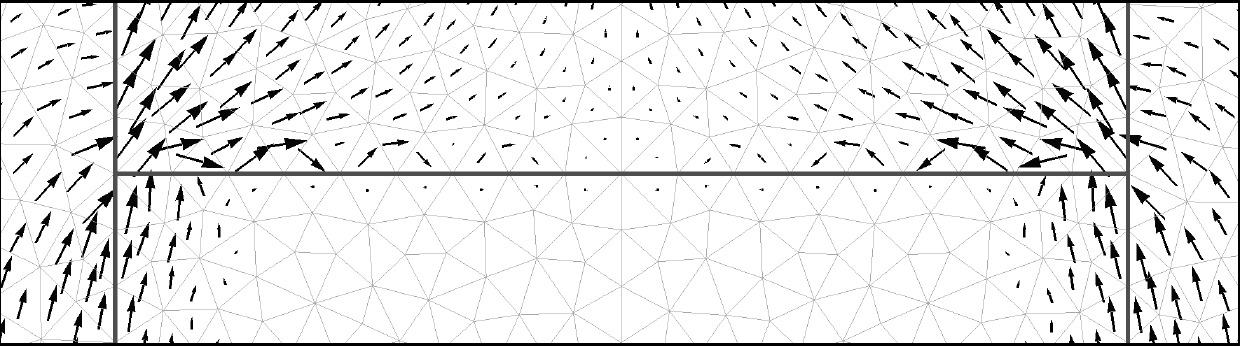}
		\caption{$\hd\in \hspdone(\Ohd)$ and $\ad\in \aspdone(\Oad)$.}
		\label{bar_zoom_a1_h1}
        \end{subfigure}
\begin{subfigure}[b]{\linewidth}  
            \centering 
		\includegraphics[width=\textwidth]{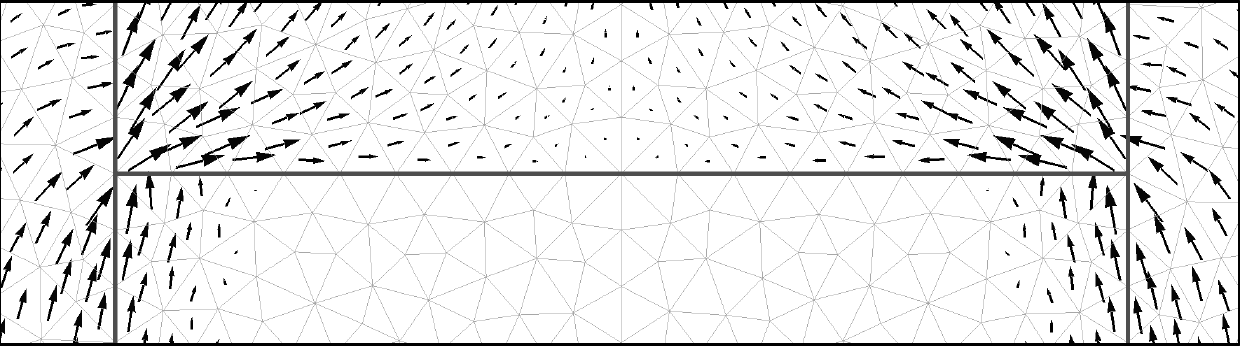}
		\caption{$\hd\in \hspdtwo(\Ohd)$ and $\ad\in \aspdone(\Oad)$.}
		\label{bar_zoom_a2_h1}	
      \end{subfigure}
        \caption{Details of two solutions for the stacked bar problem, magnetic flux density near the material interface (arrows represent the average value in each element). (a) Unstable choice of function spaces, resulting in non-physical oscillations on $\Gmd$. (b) Example of a stabilized problem with hierarchical basis functions on $\Gmd$ for $\hd$.}
        \label{bar_b_zoom}
\end{figure}

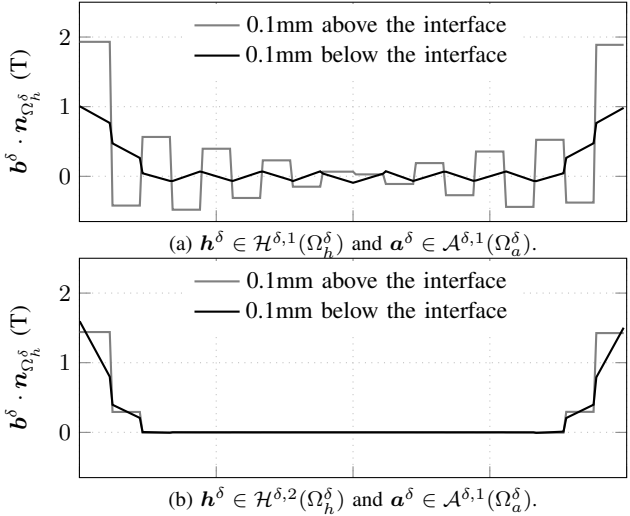
\begin{figure}[h!]
            \begin{subfigure}[b]{\linewidth}
        \centering
	\begin{tikzpicture}[trim axis left, trim axis right][font=\small]
 	\begin{axis}[
    width=1\linewidth,
    height=4.5cm,
    grid = both,
    grid style = dotted,
    xmin=-10, 
    xmax=10,
    ymin=-0.65, 
    ymax=2.5,
    xticklabel=\empty,
	xtick scale label code/.code={},
	ytick scale label code/.code={},
    ylabel={$\bd\cdot \n_{\Ohd}$ (T)},
    ylabel style={yshift=-1.5em},
    legend style={at={(0.5, 0.99)}, anchor=north, draw=none}
    ]
    \addplot[gray, thick] 
    table[x=r,y=bn_11_ferro]{stabilityMeta/bar_ha_stability.txt};
    \addplot[black, thick] 
    table[x=r,y=bn_11_super]{stabilityMeta/bar_ha_stability.txt};
    \legend{0.1mm above the interface, 0.1mm below the interface}
    \end{axis}
	\end{tikzpicture}%
	    \vspace{-0.3cm}
			\caption{$\hd\in \hspdone(\Ohd)$ and $\ad\in \aspdone(\Oad)$.}
		\label{bar_ha_stability_unstable}
	        \end{subfigure}
\begin{subfigure}[b]{\linewidth}  
        \centering
	\begin{tikzpicture}[trim axis left, trim axis right][font=\small]
 	\begin{axis}[
    width=1\linewidth,
    height=4.5cm,
    grid = both,
    grid style = dotted,
    xmin=-10, 
    xmax=10,
    ymin=-0.65, 
    ymax=2.5,
    xticklabel=\empty,
	xtick scale label code/.code={},
	ytick scale label code/.code={},
    ylabel={$\bd\cdot \n_{\Ohd}$ (T)},
    ylabel style={yshift=-1.5em},
    legend style={at={(0.5, 0.99)}, anchor=north, draw=none}
    ]
    \addplot[gray, thick] 
    table[x=r,y=bn_12_ferro]{stabilityMeta/bar_ha_stability.txt};
    \addplot[black, thick] 
    table[x=r,y=bn_12_super]{stabilityMeta/bar_ha_stability.txt};
    \legend{0.1mm above the interface, 0.1mm below the interface}
    \end{axis}
	\end{tikzpicture}%
	\vspace{-0.3cm}
		\caption{$\hd\in \hspdtwo(\Ohd)$ and $\ad\in \aspdone(\Oad)$.}
		\label{bar_ha_stability_stable}	
      \end{subfigure}
        \caption{Normal magnetic flux density distribution (horizontal position in abscissa) just above and just below the material interface for the stacked bar problem. (a) Unstable choice of function spaces, large spurious oscillations take place. (b) Stabilized solution, with higher order basis functions on $\Gmd$ for $\hd$.}
        \label{bar_ha_stability}
\end{figure}

In $h$-$a$ and $t$-$a$ coupled formulations, using lowest order Whitney elements for both fields may lead to spurious oscillations. Let us consider a typical example with the \hafOnly, in which the numerical solution exhibits non-physical oscillations when the coupling interface $\Gmd$ separates two regions of different permeabilities. The geometry is depicted in Fig.~\ref{bar_b}: two bars (height: 10 mm, width: 20 mm) are stacked and subjected to an external field. The bottom bar is a superconductor ($n=20$, $\jc=3\times 10^8$ A/m$^2$) and defines $\Ohd$, whereas the top bar is a linear ferromagnet ($\mu_\text{r} = 1000$, $\sigma_\text{ferro} = 0$). The air and ferromagnetic domains constitute $\Oad$. With $\hd\in \hspdone(\Oh)$ and $\ad\in \aspdone(\Oa)$, spurious oscillations of the magnetic flux can be clearly seen at the interface of the two materials, see Figs.~\ref{bar_zoom_a1_h1} and~\ref{bar_ha_stability_unstable}.

It is important to emphasize that oscillations are not a consequence of the nonlinearity of the equations. Indeed, if the superconductor is replaced by a linear conductor, stability issues remain, whatever the conductor resistivity value. Oscillations only appear in situations in which there is a permeability jump across the coupling boundary $\Gm$, irrespective of whether $\Gm$ is adjacent to a conducting material or not. The oscillation amplitude decreases when the permeability of the ferromagnet is lowered. Note that oscillations can therefore easily be missed in case of nonlinear ferromagnets, that saturate quickly at the large fields involved in many superconducting systems.

\begin{figure}[h!]
\centering
            \begin{subfigure}[b]{0.49\linewidth}
            \centering
		\includegraphics[width=\textwidth]{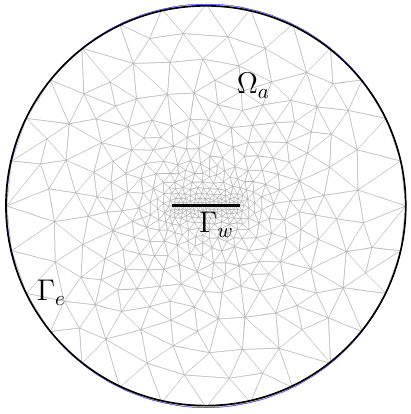}
		\caption{Problem geometry and mesh.}
		\label{tape_mesh}
        \end{subfigure}
\begin{subfigure}[b]{0.49\linewidth}  
            \centering 
		\includegraphics[width=\textwidth]{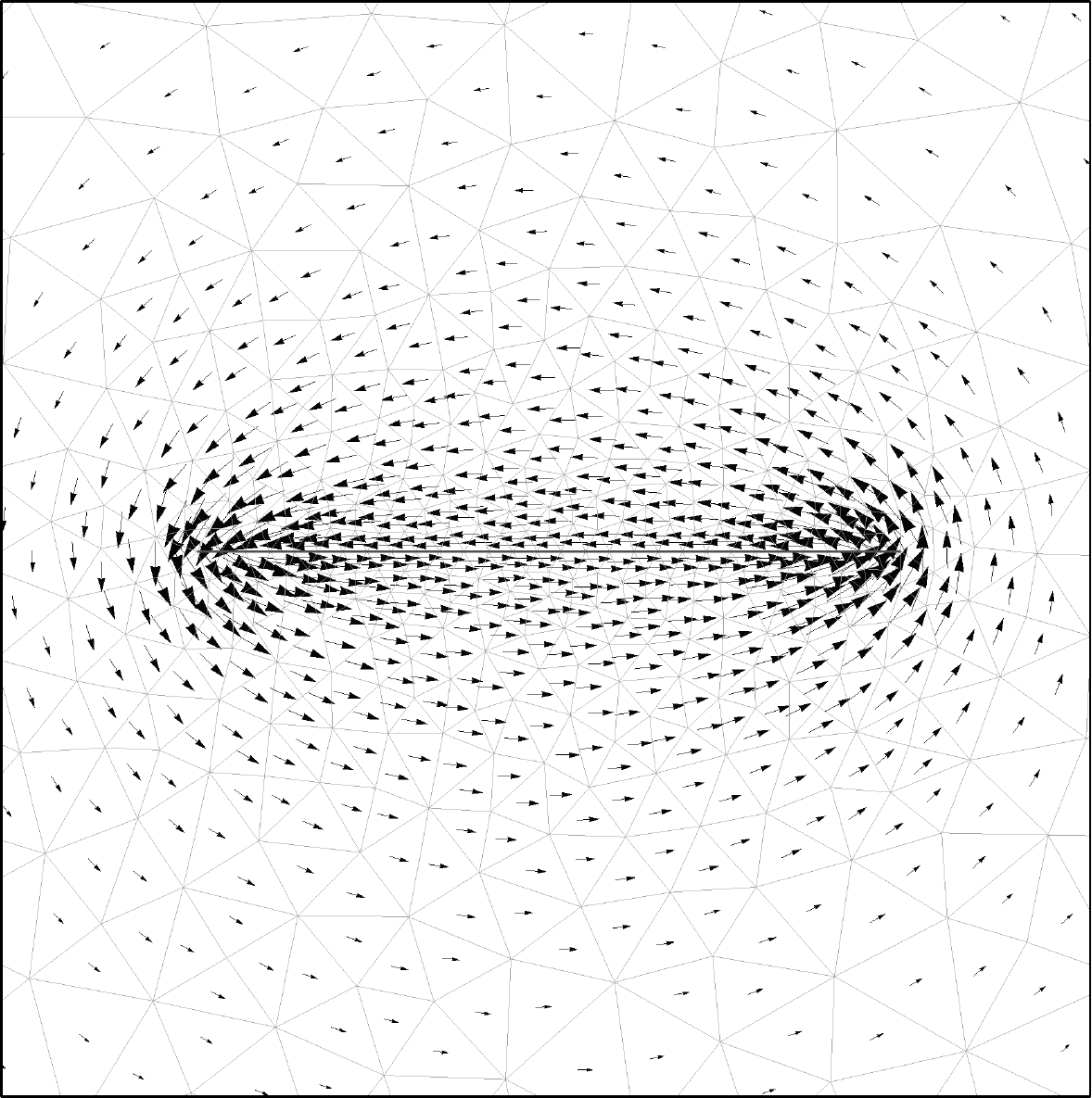}
		\caption{Magn. flux density (zoom).}
		\label{tape_zoom}	
      \end{subfigure}
        \caption{Simple problem for the \tafOnly: a superconducting tape in air, with an imposed total current intensity. (a) The problem geometry and domains. (b) Magnetic flux density in the neighbourhood of the tape, solution with first-order basis functions. Oscillations are not visible when looking at $\bd$ only.}
        \label{tape_b}
\end{figure}

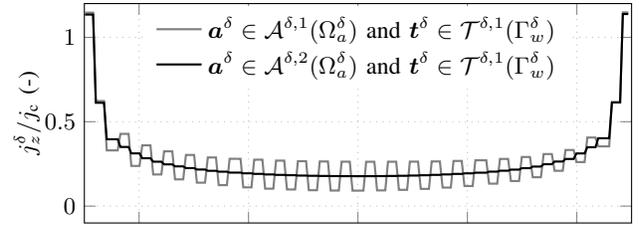
\begin{figure}[h!]
        \centering
	\begin{tikzpicture}[trim axis left, trim axis right][font=\small]
 	\begin{axis}[
    width=1\linewidth,
    height=4.5cm,
    grid = both,
    grid style = dotted,
    xmin=-5, 
    xmax=5,
    ymin=-0.1, 
    ymax=1.2,
    xticklabel=\empty,
	xtick scale label code/.code={},
	ytick scale label code/.code={},
    ylabel={$\jzd/j_\text{c}$ (-)},
    ylabel style={yshift=-1.5em},
    legend style={at={(0.5, 0.99)}, anchor=north, draw=none}
    ]
    \addplot[gray, thick] 
    table[x=r,y=j_unstable]{stabilityMeta/tape_ta_stability.txt};
    \addplot[black, thick] 
    table[x=r,y=j_stable]{stabilityMeta/tape_ta_stability.txt};
    \legend{$\ad\in\aspdone(\Oad)$ and $\td\in \tspdone(\Gwd)$, $\ad\in\aspdtwo(\Oad)$ and $\td\in \tspdone(\Gwd)$}
    \end{axis}
	\end{tikzpicture}%
        \caption{Current density for the simple tape problem. Non-physical oscillations appear when using lowest-order elements for both fields. Enriching the space for $\ad$ on $\Gwd$ stabilizes the problem and spurious oscillations disappear.}
        \label{tape_j}
\end{figure}

Similarly, in the numerical solution of the \tafOnly, the current density displays oscillations across the tape, typically at low imposed currents. A representative example is shown in Figs.~\ref{tapeGeometry_ta_article_2D} and \ref{tape_mesh}, illustrating a single straight tape in air (thickness: $10^{-3}$ mm, width: 10 mm, $n=20$, $\jc=2.5\times 10^8$ A/m$^2$), with an imposed current intensity. The magnetic flux density distribution does not exhibit problematic oscillations, but the current density is clearly non-physical. See Figs.~\ref{tape_zoom} and \ref{tape_j}. Again, such oscillations also appear with linear (ohmic) materials. The oscillation amplitude however decreases when the resistivity increases.

Oscillations can be avoided by choosing adapted function spaces. We investigate one possibility in the following subsection.

\subsection{Enriched spaces with hierarchical basis functions}

One possibility to stabilize the problem is to enrich locally the function space of one of the two fields for the \hafOnly. This is illustrated in Figs.~\ref{bar_b_zoom} and \ref{bar_ha_stability}, where using higher-order basis functions for $\hd$ on $\Gmd$ allows to overcome non-physical oscillations. Likewise, enriching the $\ad$ space yields a similar effect. For the \tafOnly, enriching the $\ad$ space stabilizes the problem, see Fig.~\ref{tape_j}. This procedure is a solution inspired by well-known results in mixed formulations in mechanics, such as Stokes' (nearly) incompressible flow problems \cite{brezzi1990discourse}. However, the situation is not the same, since fields in these problems are coupled inside the domain, whereas here, we couple the fields via boundaries of domains.

The above observations and the stability results of the next sections motivate the use of higher-order functions. Since we will only enrich functions locally, we use hierarchical functions \cite{zienkiewicz2005finite,geuzaine2001high}. In 2D, we formally associate these functions with edges. Let $n_{e,1}$ and $n_{e,2}$ be the nodes at the ends of an edge $e$, we define the associated hierarchical basis function $\psi_{2,e} \Def \psi_{n_{e,1}} \psi_{n_{e,2}}$. It vanishes on all nodes and is referred to as a bubble function.

To obtain stable formulations, it is sufficient to introduce these functions on the domain interfaces only: $\Gmd$ or $\Gwd$. We add to the expansion \eqref{eqn_h_decomposition} for $\hd$ the term $\sum_{e\in \Gmd} \phi_{2,e}\, \grad \psi_{2,e}$,
with the new degrees of freedom $\phi_{2,e}$. Note that as $\curl \grad \cdot = \vec 0$, the new terms do not contribute to the current density. We denote the resulting function space by $\hspdtwo(\Ohd)$.

For the magnetic vector potential $\ad$, we add to \eqref{eqn_a_decomposition} the term $\sum_{e\in \Gd} a_{2,e}\, \psi_{2,e} \vec z$,
with $\Gd = \Gmd$ in the \haf and $\Gd = \Gwd$ in the \tafOnly, $\vec z$ being the out-of-plane direction, and $a_{2,e}$ representing the new degrees of freedom. The associated function space is denoted by $\aspdtwo(\Oad)$.

Finally, although it will lead to issues with the Newton-Raphson technique, we consider for the current vector potential $\td$ in 2D the term $\sum_{e\in \Gamma_w} t_{2,e}\, \psi_{2,e}\n$,
in addition to decomposition \eqref{eqn_t_decomposition}, with $t_{2,e}$ being the new degrees of freedom, to construct the space $\tspdtwo(\Gwd)$.

As explained in what follows, these hierarchical basis functions enrich the polynomial order of the span of function traces on $\Gmd$ and $\Gwd$ by one. This will be shown to be sufficient to obtain stability in sections \ref{sec_haf} and \ref{sec_taf}.

\subsection{Function space for the traces on $\Gmd$ and $\Gwd$}

The coupling integrals involve the traces of functions on the interfaces $\Gmd$ and $\Gwd$. The range of these traces depends on the chosen function spaces and, as we will see, determines the stability of the system.

The trace $\hd\times \n|_{\Gmd}$ of $\hd\in \hspdone(\Ohd)$, involved in the coupling terms of the \haf is in the $\vec z$-direction and is piecewise constant. With the hierarchical enrichment, it becomes piecewise linear. Similarly, in the \tafOnly, with $\td\in \tspdone(\Gwd)$, $\curl \t$ is along $\vec z$ and is piecewise constant. For $\td\in \tspdtwo(\Gwd)$, $\curl\t$ is piecewise linear.

The vector potential $\ad \in\aspdone(\Oad)$ is along $\vec z$. It is continuous and piecewise linear. For $\a \in\aspdtwo(\Oad)$, it is continuous and piecewise quadratic.

Now that the discrete framework has been presented, we summarize the basics of saddle point stability analysis, and then apply the theory on the two-bar and single-tape examples, in order to explain how a proper choice of approximation space leads to a stable coupled formulation.

\section{Basics of Stability Analysis}\label{sec_stabilityAnalysis}

Mixed finite element formulations face numerical stability issues if function spaces for their unknowns are not chosen consistently. Typically, spurious oscillations in the numerical solution may arise and affect the accuracy of the method, as illustrated in the previous section. The theory of mixed finite element formulations \cite{brezziBook} provides compatibility conditions on spaces to ensure the numerical stability of the problem.

Both the \haf and \taf fit into the classical framework of perturbed saddle-point problems, as will be shown in sections \ref{sec_haf} and \ref{sec_taf}. In this section, we present the stability conditions relevant to this class of problems, following closely Brezzi's classical treatment \cite{brezziBook}. We also describe the \textit{inf-sup test} \cite{chapelle1993inf}, a numerical test which is used for checking the compatibility of specific functions spaces in the discrete setting.

\subsection{Theory}

We consider two Hilbert spaces $V$ and $Q$ and their dual space $V'$ and $Q'$ containing all linear functionals $V\to \mathbb{R}$ and $Q\to \mathbb{R}$, respectively. On these spaces, we build perturbed saddle-point problems of the following form: for given $f\in V'$ and $g\in Q'$, find $u\in V$ and $p\in Q$ such that
\begin{equation}\label{eqn_perturbedProblem}
\left\{
\begin{aligned}
a(u,v) + b(v,p) &= \langle f,v\rangle,\quad \forall v\in V,\\
b(u,q) - c(p,q) &= \langle g,q\rangle,\quad \forall q\in Q,
\end{aligned}\right.
\end{equation}
with $a(\cdot,\cdot)$, $b(\cdot,\cdot)$, and $c(\cdot,\cdot)$ continuous bilinear forms on $V\times V$, $V \times Q$, and $Q \times Q$, respectively, and where $\langle f, v\rangle$ (resp. $\langle g, q\rangle$) denotes the value of the functional $f$ (resp. $g$) at $v$ (resp. $q$). The term $-c(p,q)$ is considered to be a perturbation of the classical saddle-point problem.

In practice, we solve a discretized version of \eqref{eqn_perturbedProblem} and look for $u^\delta$ and $p^\delta$ in finite-dimensional spaces $V^\delta\subseteq V$ and $Q^\delta\subseteq Q$ respectively (the finite element spaces), with operators defined on these discrete spaces.

Spaces $V^\delta$ and $Q^\delta$ are equipped with norms $\|\cdot \|_{V^\delta}$ and $\|\cdot \|_{Q^\delta}$. Dual norms are used for elements in the dual spaces ${V^{\delta}}'$ and ${Q^{\delta}}'$. The norms of the bilinear operators are finite and are defined as follows:
\begin{align}
\|a^\delta\| &\Def\sup_{u^\delta,v^\delta\in V^\delta} \frac{a(u^\delta,v^\delta)}{\|u^\delta\|_{V^\delta}\|v^\delta\|_{V^\delta}}, \label{eqn_anorm}\\
\|b^\delta\| &\Def \sup_{u^\delta\in V^\delta, q^\delta\in Q^\delta} \frac{b(u^\delta,q^\delta)}{\|u^\delta\|_{V^\delta}\|q^\delta\|_{Q^\delta}}, \label{eqn_bnorm}\\
\|c^\delta\| &\Def \sup_{p^\delta,q^\delta\in Q^\delta} \frac{c(p^\delta,q^\delta)}{\|p^\delta\|_{Q^\delta}\|q^\delta\|_{Q^\delta}}. \label{eqn_cnorm}
\end{align}



Before stating the main stability theorem, we finally introduce the kernels
\begin{align}
K^\delta &\Def \{v^\delta\in V^\delta : b(v^\delta,q^\delta) = 0, \forall q^\delta\in Q^\delta\},\\
H^\delta &\Def \{q^\delta\in Q^\delta : b(v^\delta,q^\delta) = 0, \forall v^\delta\in V^\delta\}.
\end{align}

Let $(\bar u,\bar p)$ be the exact solution of the perturbed saddle point-problem \eqref{eqn_perturbedProblem}. With $a(\cdot,\cdot)$ and $c(\cdot,\cdot)$ symmetric, positive semi-definite continuous bilinear forms, Proposition 5.5.2. in \cite{brezziBook} states that, if $a(\cdot,\cdot)$ is coercive on $K^\delta$ and $c(\cdot,\cdot)$ is coercive on $H^\delta$, i.e., if there exists two constants $\alpha^\delta>0$ and $\gamma^\delta>0$ such that
\begin{align}
a(v^\delta,v^\delta) \ge \alpha^\delta \|v\|_{V^\delta}^2 ,\quad &\forall v^\delta\in K^\delta,\label{coercivenessoOnKernels_a}\\
c(q^\delta,q^\delta) \ge \gamma^\delta \|q\|_{Q^\delta}^2, \quad &\forall q^\delta\in H^\delta,\label{coercivenessoOnKernels_c}
\end{align}
and if there exists a constant $\beta^\delta>0$ such that
\begin{align}\label{infsupCondition}
\inf_{q^\delta\in (H^\delta)^\perp} \sup_{v^\delta\in (K^\delta)^\perp} \frac{b(v^\delta,q^\delta)}{\|q^\delta\|_{Q^\delta} \|v^\delta\|_{V^\delta}} = \beta^\delta > 0
\end{align}
is satisfied on the orthogonal complements $\cdot^\perp$ of $H^\delta$ and $K^\delta$, then the discretized problem has a unique solution $(u^\delta,p^\delta)$ satisfying
\begin{multline}\label{infsupCondition_convergence}
\|u^\delta - \bar u\|_{V^\delta} + \|p^\delta - \bar p\|_{Q^\delta} \\
\le C^\delta \paren{\inf_{v^\delta\in V^\delta} \|v^\delta-\bar u\|_{V^\delta} + \inf_{q^\delta\in Q^\delta} \|q^\delta-\bar p\|_{Q^\delta}},
\end{multline}
with a constant $C^\delta$ depending only on the stability constants, $\alpha^\delta$, $\beta^\delta$, $\gamma^\delta$, and continuity constants, $\|a^\delta\|$, $\|b^\delta\|$, $\|c^\delta\|$.

If these values can be chosen independent of the mesh, $C^\delta$ is bounded with mesh refinement and the problem is said to be stable.

Condition \eqref{infsupCondition} is the so-called inf-sup condition, or the Babu{\v{s}}ka-Brezzi condition \cite{babuvska1973finite,brezzi1974existence}.


\subsection{Numerical Inf-Sup Test}

In most practical cases, the inf-sup value $\beta^\delta$ cannot be evaluated analytically. Instead, it can be estimated with a numerical inf-sup test \cite{brezziBook,chapelle1993inf}. On a given mesh, unknown fields $v^\delta \in V^\delta$ and $q^\delta \in Q^\delta$ are described by vectors $\vec v^{\delta}$ and $\vec q^{\delta}$ containing the degrees of freedom. We introduce orthogonal matrices $\mat N^{\delta}_{V^\delta}$ and $\mat N^{\delta}_{Q^\delta}$ such that $\|v^{\delta}\|_{V^{\delta}}^2 = \paren{\vec v^{\delta}}\transpose \mat N^{\delta}_{V^\delta} \vec v^{\delta}$ and $\|q^{\delta}\|_{Q^{\delta}}^2 = \paren{\vec q^{\delta}}\transpose \mat N^{\delta}_{Q^{\delta}}\vec q^{\delta}$ and we introduce $\mat B^{\delta}$, the coupling matrix satisfying $b(v^{\delta},q^{\delta}) = \paren{\vec q^{\delta}}\transpose \mat B^{\delta} \vec v^{\delta}$, obtained from the finite element assembly.

In terms of these matrices, Eq.~\eqref{infsupCondition} reads
\begin{align}\label{infsupCondition_matrix}
\inf_{\vec q^{\delta}\in (H^\delta)^\perp} \sup_{\vec v^{\delta}\in (K^\delta)^\perp} \frac{\paren{\vec q^{\delta}}\transpose \mat B^{\delta} \vec v^{\delta}}{ \paren{\paren{\vec q^{\delta}}\transpose \mat N^{\delta}_{Q^{\delta}}\vec q^{\delta}} \paren{\paren{\vec v^{\delta}}\transpose \mat N^{\delta}_{V^\delta} \vec v^{\delta}}} = \beta^\delta.
\end{align}
The inf-sup value $\beta^{\delta}$ in Eq.~\eqref{infsupCondition_matrix} can be shown to be equal to the square root of the smallest non-zero eigenvalue of the generalized eigenvalue problem \cite{malkus1981eigenproblems}
\begin{align}\label{eqn_eigenvalueProblem}
\paren{\mat B^{\delta} \paren{\mat N^{\delta}_{V^{\delta}}}^{-1} \paren{\mat B^{\delta}}\transpose} \vec q^{\delta} &= \lambda^{\delta} \mat N^{\delta}_{Q^{\delta}} \vec q^{\delta}.
\end{align}
Note that we disregard zero eigenvalues because they are associated with eigenvectors defining elements in $H^\delta$ that are not involved in the inf-sup condition. Note also that the norm $\|b^\delta\|$ is the square root of the largest eigenvalue of problem \eqref{eqn_eigenvalueProblem}.

The inf-sup test consists in computing $\beta^\delta$ values for progressively refined meshes. If the values appear to be bounded from below by a positive value independent of mesh size, and if the other conditions (coerciveness and continuity) are met, then the sequence of problems is considered to be stable. On the other hand, if some eigenvalues tend to zero, we expect stability issues, because the inf-sup condition then fails to be satisfied. Even if the numerical test does not provide a formal proof of stability, experience shows that it is a reliable indicator \cite{chapelle1993inf,bathe2001inf}.

In the next sections we apply these results on the $h$-$a$- and $t$-$a$-formulations in order to investigate the stability.

From now on, we only keep the $\cdot^\delta$ superscript for functions spaces, continuity, coerciveness and inf-sup values, to stress the importance of mesh-dependency. As we will only stay in the discrete setting, we drop it elsewhere, for conciseness.

\section{Analysis of the \haf}\label{sec_haf}

For simplicity, we start by presenting the \haf on materials with constant permeability and conductivity. We will then extend the conclusions to systems with superconductors and nonlinear ferromagnetic materials.

The analysis is restricted to 2D problems with in-plane magnetic field.

\subsection{Linear materials}

We begin the analysis with a linear problem, i.e., materials that have a constant resistivity and reluctivity, but are not necessarily homogeneous.

Using the implicit Euler method, at a given time step $n$, the solution $(\a,\h) \Def (\a_n, \h_n)$ depends on the solution at the previous time step $(\cdot)_{n-1}$. If we multiply the first equation of \eqref{eqn_coupledFormulation} by the time step $\Delta t$, we obtain the system
\begin{equation}\label{eqn_coupledFormulation_lin_oneSolve}
\begin{aligned}
& \volIntBig{\mu\, \h}{\h'}{\Oh} + \volInt{\Delta t\, \rho\, \curl \h}{\curl \h'}{\Ohc}\\
& \qquad\qquad + \surInt{\a \times\vec n_{\Oh}}{\h'}{\Gamma_\text{m}} = \langle \vec s, \h' \rangle,\\
&\surInt{\h\times\vec n_{\Oa}}{\a'}{\Gamma_\text{m}} - \volInt{\nu\, \curl \a}{\curl \a'}{\Oa}= 0,
\end{aligned}
\end{equation}
with the right-hand side functional defined by
\begin{multline}
\langle \vec s, \h' \rangle = \surInt{\a_{n-1} \times\vec n_{\Oh}}{\h'}{\Gamma_\text{m}} + \volIntBig{(\mu\, \h)_{n-1}}{\h'}{\Oh}\\
 - \Delta t \sum_{i\in C} V_i \mathcal{I}_i(\h').
\end{multline}
System \eqref{eqn_coupledFormulation_lin_oneSolve} can be rewritten as
\begin{equation}\label{ha_LinearIdentificationWithTheory}
\begin{aligned}
& \volIntBig{\mu\ \h}{\h'}{\Oh} + \volInt{\Delta t\ \rho\ \curl \h}{\curl \h'}{\Ohc}\\
& \qquad  + \surInt{\a \times\vec n_{\Oh}}{\h'}{\Gamma_\text{m}} = \langle \vec s, \h' \rangle,\\
& \surInt{\a'\times\vec n_{\Oh}}{\h}{\Gamma_\text{m}} - \volInt{\nu\ \curl \a}{\curl \a'}{\Oa} = 0,
\end{aligned}
\end{equation}
using $\vec n_{\Oa} = - \vec n_{\Oh}$. For conciseness, we consider homogeneous essential boundary conditions. Problem \eqref{ha_LinearIdentificationWithTheory} can be cast into the form of Eq.~\eqref{eqn_perturbedProblem}, with identical function spaces for unknown functions and test functions. The case of non-homogeneous essential boundary conditions can be easily treated, and the analysis remains unchanged.

After discretization, we obtain a system of linear equations in a matrix-vector form. The formulation will be considered stable if a sequence of problems on progressively refined meshes satisfies conditions \eqref{eqn_anorm} to \eqref{eqn_cnorm} and \eqref{coercivenessoOnKernels_a} to \eqref{infsupCondition}, with constants $\alpha^\delta$, $\beta^\delta$, $\gamma^\delta$, $\|a^\delta\|$, $\|b^\delta\|$, and $\|c^\delta\|$ independent of mesh size.

In $\hspzd(\Oh)$ and $\aspzd(\Oa)$, we define the norms
\begin{align}\label{ha_normH}
\|\h\|_{\hspzd}^2 &= \volInt{\mu_0\,\h}{\h}{\Oh}+ \volInt{\Delta t_0\,\rho_0\,\curl \h}{\curl \h}{\Ohc},\\
\|\a\|_{\aspzd}^2 &= \volInt{\nu_0\,\curl \a}{\curl \a}{\Oa},\label{ha_normA}
\end{align}
with $\rho_0$ being a characteristic resistivity (e.g., the resistivity of region $\Ohc$) and $\Delta t_0$ a characteristic time step. With these norms, whatever the discretization, $\forall \h\in \hspzd(\Oh)$ and $\forall \a\in \aspzd(\Oa)$,
\begin{align}
a(\h,\h) &= \volIntBig{\mu\, \h}{\h}{\Oh} + \volInt{\Delta t\, \rho\, \curl \h}{\curl \h}{\Ohc}\notag\\
&\ge \min(\mu/\mu_0,\Delta t/\Delta t_0\cdot \rho/\rho_0)\, \|\h\|_{\hspzd}^2,\\
c(\a,\a) &= \volInt{\nu\, \curl \a}{\curl \a}{\Oa}\notag\\
&\ge \min(\nu/\nu_0)\, \|\a\|_{\aspzd}^2,
\end{align}
which proves the coerciveness properties \eqref{coercivenessoOnKernels_a}, and \eqref{coercivenessoOnKernels_c}, with $\alpha^\delta \ge \min(\mu/\mu_0,\Delta t/\Delta t_0\cdot \rho/\rho_0)>0$, and $\gamma^\delta \ge \min(\nu/\nu_0)>0$. Similarly, we can prove $\|a^\delta\| \le \max(\mu/\mu_0,\Delta t/\Delta t_0\cdot \rho/\rho_0)<\infty$, and $\|c^\delta\| \le \max(\nu/\nu_0)<\infty$, using the Cauchy-Schwarz inequality.

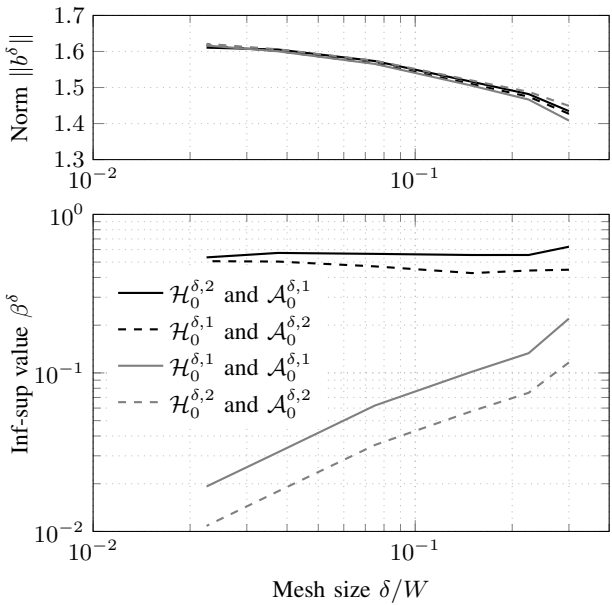
\begin{figure}[h!]
  \begin{subfigure}[b]{0.49\textwidth}  
        \centering
	\begin{tikzpicture}[trim axis left, trim axis right][font=\small]
 	\begin{semilogxaxis}[
	tick scale binop=\times,
    width=0.95\linewidth,
    height=3.5cm,
    grid = both,
    grid style = dotted,
    xmin=0.01, 
    xmax=0.4,
    ymin=1.3, 
    ymax=1.7,
	xtick scale label code/.code={},
    xtick={0.01, 0.1, 1},
	ytick scale label code/.code={},
    ylabel={Norm $\|b^\delta\|$},
    ylabel style={yshift=-0.8em},
    legend style={at={(0.03, 0.35)}, anchor=west, draw=none}
    ]
    \addplot[black, thick] 
    table[x=meshsize,y=normb12]{stabilityMeta/infsup_ha_stackedBarLinear.txt};
        \addplot[black, dashed, thick] 
    table[x=meshsize,y=normb21]{stabilityMeta/infsup_ha_stackedBarLinear.txt};
    \addplot[gray, thick] 
    table[x=meshsize,y=normb11]{stabilityMeta/infsup_ha_stackedBarLinear.txt};
        \addplot[gray, dashed, thick] 
    table[x=meshsize,y=normb22]{stabilityMeta/infsup_ha_stackedBarLinear.txt};
    \end{semilogxaxis}
\end{tikzpicture}%
\end{subfigure}
        \hfill
\begin{subfigure}[b]{0.49\textwidth}
        \centering
	\begin{tikzpicture}[trim axis left, trim axis right][font=\small]
 	\begin{loglogaxis}[
	tick scale binop=\times,
    width=0.95\linewidth,
    height=5.8cm,
    grid = both,
    grid style = dotted,
    xmin=0.01, 
    xmax=0.4,
    ymin=0.01, 
    ymax=1,
	xtick scale label code/.code={},
    xtick={0.01, 0.1, 1},
	ytick scale label code/.code={},
	xlabel={Mesh size $\delta/W$},
    ylabel={Inf-sup value $\beta^\delta$},
    ylabel style={yshift=-0.8em},
    legend style={at={(0.03, 0.58)}, anchor=west, draw=none}
    ]
    \addplot[black, thick] 
    table[x=meshsize,y=beta12]{stabilityMeta/infsup_ha_stackedBarLinear.txt};
        \addplot[black, dashed, thick] 
    table[x=meshsize,y=beta21]{stabilityMeta/infsup_ha_stackedBarLinear.txt};
    \addplot[gray, thick] 
    table[x=meshsize,y=beta11]{stabilityMeta/infsup_ha_stackedBarLinear.txt};
        \addplot[gray, dashed, thick] 
    table[x=meshsize,y=beta22]{stabilityMeta/infsup_ha_stackedBarLinear.txt};
    \legend{$\hspzdtwo$ and $\aspzdone$,$\hspzdone$ and $\aspzdtwo$,$\hspzdone$ and $\aspzdone$,$\hspzdtwo$ and $\aspzdtwo$}
    \end{loglogaxis}
	\end{tikzpicture}%
	  \end{subfigure}
        \caption{Evolution of the inf-sup constant $\beta^\delta$ from Eq.~\eqref{eqn_infsup_ha} and norm $\|b\|$ with mesh refinement ($\delta\to 0$) on the stacked bar linear problem. Four cases are considered: $\h \in \hspzdi(\Oh)$ and $\a \in \aspzdj(\Oa)$, for $(i,j)\in \{1,2\}\times \{1,2\}$. We can only conclude on stability when $i\neq j$, i.e., when exactly one space is enriched with respect to Whitney elements (black lines).}
        \label{infsup_ha_stackedBarLinear}
\end{figure}

To guarantee stability, the inf-sup condition remains to be met. There must exist a $\beta^\delta > 0$ independent of mesh size that fulfils
\begin{align}\label{eqn_infsup_ha}
\inf_{\a \in H^{\perp}} \sup_{\h \in K^\perp} \frac{\surInt{\a \times\vec n_{\Oh}}{\h}{\Gamma_\text{m}}}{\|\a\|_{\aspzd} \|\h\|_{\hspzd}} \ge \beta^\delta.
\end{align}
We also have to verify that $\|b^\delta\|$ is bounded from above. To check both properties, a numerical inf-sup test is conducted on the stacked bar geometry represented in Fig.~\ref{bar_b}, with linear homogeneous materials ($\rho = 1.6\times 10^{-8}$ $\Omega$m, $\mur = 1000$, non-conducting ferromagnet), for different discretization choices. Results are shown in Fig.~\ref{infsup_ha_stackedBarLinear}, with $\rho_0 = 1.6\times 10^{-8}$~$\O$m.

First, the norm $\| b^\delta\|$ of the coupling operator is bounded from above independent of the function spaces, as shown in the upper part of Fig.~\ref{infsup_ha_stackedBarLinear}. However, the evolution of the inf-sup value shows two different behaviors. When exactly one of the two fields $\h$ and $\a$ is enriched with hierarchical elements, the inf-sup value does not decrease and \eqref{infsupCondition_convergence} ensures the stability of the associated problem. Otherwise, the inf-sup value typically decreases with $\beta^\delta \sim \delta$ and stability issues, i.e., oscillations in the numerical solution, are expected.

In practice, when choosing $\h \in \hspzdone(\Oh)$ and $\a \in \aspzdone(\Oa)$, or $\h \in \hspzdtwo(\Oh)$ and $\a \in \aspzdtwo(\Oa)$, we do observe such oscillations. However, they only appear at interfaces with large permeability jumps. In contrast, when no ferromagnetic material is present in the geometry, the numerical results are satisfying. These behaviors can be explained by Proposition 4.3.1 of \cite{brezziBook}, that follows from the Lax-Milgram theorem. When $a(\cdot,\cdot)$ and $c(\cdot,\cdot)$ are coercive, irrespective of whether the coupling operator satisfies the inf-sup condition, the problem has a unique solution $(\a,\h)$ and we have the following inequality:
\begin{align}\label{eqn_stabilityPerturbed}
\frac{\alpha^\delta}{2}\|\a \|_{\aspzd} + \frac{\gamma^\delta}{2}\| \h \|_{\hspzd} \le \frac{1}{2\alpha^\delta} \|\vec{s}_{\a} \|_{(\aspzd)'} +\frac{1}{2\gamma^\delta} \| \vec{s}_{\h} \|_{(\hspzd)'},
\end{align}
with $\vec{s}_{\a}$ and $\vec{s}_{\h}$ the right-hand sides of the final system (after treating non-homogeneous essential boundary conditions). The problem is actually always stable in the sense of Eq.~\eqref{eqn_stabilityPerturbed}. However, the provided bound deteriorates when either $\alpha^\delta$ or $\gamma^\delta$ decreases, which is the case when considering a ferromagnetic material in $\Oa$. Indeed, when $1/\nu \to \infty$ in $\Oa$, $\gamma^\delta\to 0$. With practical mesh resolutions, the bound in Eq.~\eqref{eqn_stabilityPerturbed} is not strict enough and stability issues arise.

As said above, we can extend to spaces with non-homogeneous essential boundary conditions. In practice, when dealing with ferromagnetic materials adjacent to $\Gm$, it is therefore recommended to choose either $\h \in \hspdone(\Oh)$ and $\a \in \aspdtwo(\Oa)$, or $\h \in \hspdtwo(\Oh)$ and $\a \in \aspdone(\Oa)$, to guarantee stability.

\subsection{Interpretation}

To illustrate the link between the eigenvalue problem and the spurious oscillations, we investigate Eq.~\eqref{eqn_eigenvalueProblem} applied on the stacked-bar problem with linear materials.

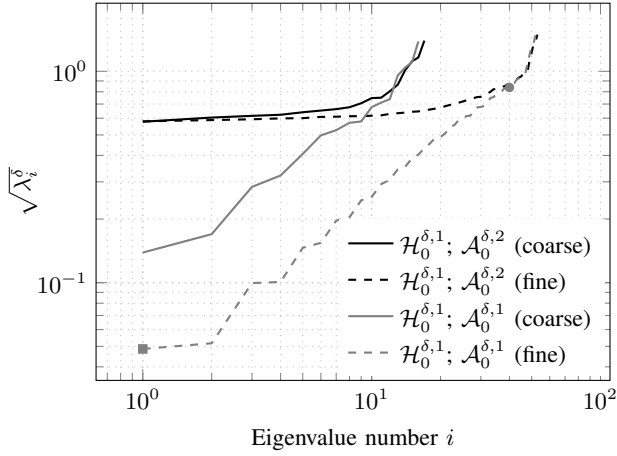
\begin{figure}[h!]
        \centering
	\begin{tikzpicture}[trim axis left, trim axis right][font=\small]
 	\begin{loglogaxis}[
	tick scale binop=\times,
    width=0.95\linewidth,
    height=6.6cm,
    grid = both,
    grid style = dotted,
    xmax=110,
	xtick scale label code/.code={},
	ytick scale label code/.code={},
	xlabel={Eigenvalue number $i$},
    ylabel={$\sqrt{\lambda^{\delta}_i}$},
    ylabel style={yshift=-0.8em},
    legend style={at={(0.99, 0.01)}, cells={anchor=west}, anchor=south east, draw=none},
    ]
    \addplot[black, thick] 
    table[x=ID,y=sqrt_eigenvalue]{stabilityMeta/infsup_ha_stackedBar_sqrt_eigenvalues_12_coarser.txt};
            \addplot[black, dashed, thick] 
    table[x=ID,y=sqrt_eigenvalue]{stabilityMeta/infsup_ha_stackedBar_sqrt_eigenvalues_12_coarse.txt};
    \addplot[gray, thick] 
    table[x=ID,y=sqrt_eigenvalue]{stabilityMeta/infsup_ha_stackedBar_sqrt_eigenvalues_11_coarser.txt};
    \addplot[gray, dashed, thick] 
    table[x=ID,y=sqrt_eigenvalue]{stabilityMeta/infsup_ha_stackedBar_sqrt_eigenvalues_11_coarse.txt};
    \addplot[only marks, gray, mark=square*, mark options={gray, scale=0.8, style={solid}}] 
    coordinates {(1, 0.0486)};
        \addplot[only marks, gray, mark=*, mark options={gray, scale=0.8, style={solid}}] 
    coordinates {(40, 0.8382)};
    \legend{$\hspzdone$; $\aspzdtwo$ (coarse),$\hspzdone$; $\aspzdtwo$ (fine), $\hspzdone$; $\aspzdone$ (coarse),$\hspzdone$; $\aspzdone$ (fine)}
    \end{loglogaxis}
	\end{tikzpicture}%
	        \caption{Distribution of the square root of the non-zero eigenvalues from problem \eqref{eqn_eigenvalueProblem} on the stacked-bar geometry with linear materials. The smallest values are the inf-sup values $\beta^\delta$, the largest are the norms $\|b^\delta\|$. Eigenvectors associated with the square and circle points are represented in Fig.~\ref{bar_b_zoom_eigen}.}
        \label{infsup_ha_eigenvalues}
\end{figure}

\begin{figure}[h!]
\centering
            \begin{subfigure}[b]{\linewidth}
            \centering
		\includegraphics[width=0.85\textwidth]{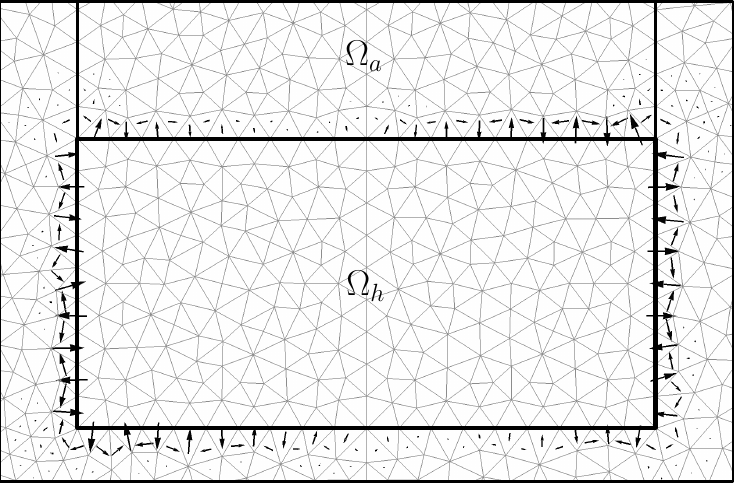}
		\caption{Associated with the square in Fig.~\ref{infsup_ha_eigenvalues}, smallest non-zero eigenvalue.}
		\label{bar_zoom_a1_h1_eigen_min}
        \end{subfigure}
\begin{subfigure}[b]{\linewidth}  
            \centering 
		\includegraphics[width=0.85\textwidth]{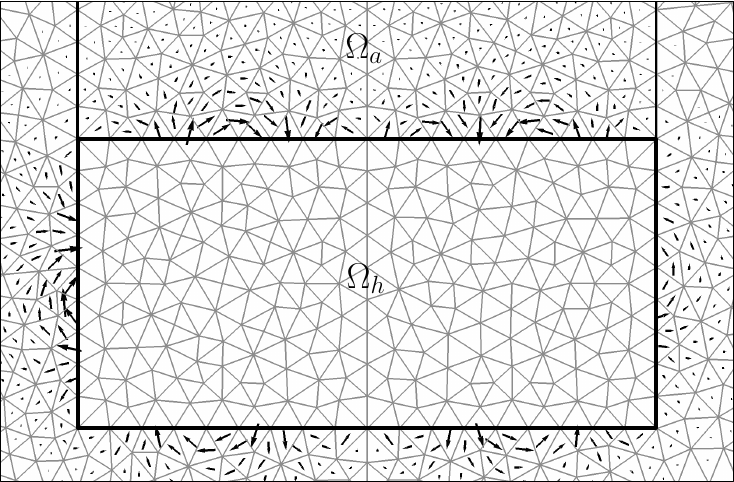}
		\caption{Associated with the circle in Fig.~\ref{infsup_ha_eigenvalues}.}
		\label{bar_zoom_a1_h1_eigen_40}	
      \end{subfigure}
        \caption{Eigenvectors associated with the two dots in Fig.~\ref{infsup_ha_eigenvalues}, for $\h\in \hspdone(\Oh)$ and $\a\in \aspdone(\Oa)$ (unstable) and the same mesh as in Fig.~\ref{bar_b_zoom}. The thick curve is $\Gm$.}
        \label{bar_b_zoom_eigen}
\end{figure}

For the stability analysis, we are only interested in the non-zero eigenvalues, whose associated eigenvectors form a basis of $H^\perp$. These eigenvalues are represented in Fig.~\ref{infsup_ha_eigenvalues} for both stable and unstable choices of function spaces, at two different discretization levels. The same conclusions as from Fig.~\ref{infsup_ha_stackedBarLinear} can be drawn, by looking only at the smallest eigenvalue. In particular, with the choice $\h\in \hspdone(\Oh)$ and $\a\in \aspdone(\Oa)$, the problem is unstable because it contains modes of smaller and smaller eigenvalues when the mesh is refined. The eigenvector associated with the smallest eigenvalue is represented in Fig.~\ref{bar_zoom_a1_h1_eigen_min}. Clearly, such a mode (among others) is also activated in the unstable solution of Fig.~\ref{bar_zoom_a1_h1} with nonlinear materials. Its weight in the coupling term $\surInt{\a\times\n_{\Oh}}{\h}{\Gm}$ is small with respect to its norm.

Of course, such oscillating modes still exist in the $H^\perp$ basis with the stable choice $\h\in \hspdone(\Oh)$ and $\a\in \aspdtwo(\Oa)$, but their eigenvalues have been leveled up and new modes not longer introduce smaller and smaller eigenvalues.

\subsection{Nonlinear materials}

We consider anhysteretic ferromagnets in $\Oa$, characterized by a saturation law for the permeability, and type-II superconductors in $\Oh$, whose resistivity is described by a power law. The associated system of equations after time discretization is as in Eq.~\ref{eqn_coupledFormulation_lin_oneSolve}, but with variable coefficients $\rho$ and $\nu$. With a Newton-Raphson linearization, we obtain a problem that is iteratively solved. The solution $(\h,\a)\Def (\h_n^k,\a_n^k)$ at time step $n$ and iteration $k$ depends on the solutions at the previous time step $(\cdot)_{n-1}$ and previous iteration $(\cdot)^{k-1}$. Using, $\vec n_{\Oa} = - \vec n_{\Oh}$, we obtain the linear system
\begin{equation}\label{oneIterationHASystem}
\begin{aligned}
& \volIntBig{\mu_0\, \h}{\h'}{\Oh} + \volInt{\Delta t\paren{\partial\e/\partial\j}^{k-1}\  \curl \h}{\curl \h'}{\Ohc}\\
& \qquad  + \surInt{\a \times\vec n_{\Oh}}{\h'}{\Gamma_\text{m}} = \langle \tilde{\vec s}_{\h}, \h' \rangle,\\
& \surInt{\a'\times\vec n_{\Oh}}{\h}{\Gamma_\text{m}} - \volInt{(\partial \h/\partial \b)^{k-1}\, \curl \a}{\curl \a'}{\Oa}\\
& \qquad \qquad \qquad\qquad\qquad = \langle \tilde{\vec s}_{\a}, \a' \rangle,
\end{aligned}
\end{equation}
with right-hand side functionals $\tilde{\vec s}_{\h}$ and $\tilde{\vec s}_{\a}$ defined by
\begin{align}
\begin{aligned}
&\langle \tilde{\vec s}_{\h}, \h' \rangle = \surInt{\a_{n-1} \times\vec n_{\Oh}}{\h'}{\Gamma_\text{m}} + \volIntBig{(\mu\, \h)_{n-1}}{\h'}{\Oh}\\
&\qquad  - \volInt{\Delta t \paren{\paren{\rho \mat I - \partial\e/\partial\j} \curl \h}^{k-1}}{\curl \h'}{\Ohc}\\
&\qquad - \Delta t \sum_{i\in C} V_i \mathcal{I}_i(\h'),
\end{aligned}\\
\langle \tilde{\vec s}_{\a}, \a' \rangle = - \volInt{ \paren{\paren{\nu \mat I- \partial\h/\partial\b} \curl \a}^{k-1}}{\curl \a'}{\Oa}
\end{align}
with the identity matrix $\mat I$. The structure is similar to that of system \eqref{ha_LinearIdentificationWithTheory}. Coerciveness and continuity of diagonal operators, $a(\cdot,\cdot)$ and $c(\cdot,\cdot)$, in the sense of norms \eqref{ha_normH} and \eqref{ha_normA} are only satisfied if the eigenvalues of matrices $\paren{\partial\e/\partial\j}^{k-1}$ and $\paren{\partial\h/\partial\b}^{k-1}$ are bounded away from zero and infinity, independently of the mesh. This is the case for the differential reluctivity with classical saturation laws. However, using the power law, the differential resistivity tends to zero for small current densities so that we cannot verify the coerciveness condition with norm \eqref{ha_normH}. Note that continuity is not satisfied either.

As for the inf-sup value $\beta^\delta$ and norm $\|b^\delta\|$, results are exactly similar to those in Fig.~\ref{infsup_ha_stackedBarLinear}. In contrast to the linear case, we do not establish a formal proof of stability due to this particular operator $a(\cdot,\cdot)$. However, we found that in practice, the conclusions obtained for the linear case remain and lead to the same recommendations. When choosing $\h \in \hspdone(\Oh)$ and $\a \in \aspdtwo(\Oa)$, or $\h \in \hspdtwo(\Oh)$ and $\a \in \aspdone(\Oa)$, we observe stable results, whereas the other combinations lead to spurious oscillations. Note that in the large fields involved with high-temperature superconductors, the ferromagnets usually saturate quickly, and the oscillation amplitude decreases.

To avoid the technical difficulty due to the power law, we could use a regularized version, 
with two limiting resistivity values. 
See also \cite{van2015numerical,laforest2018p} for a rigorous treatment of the power law in simpler formulations.

\section{Analysis of the \taf}\label{sec_taf}

We directly consider a nonlinear material in $\Gamma_w$, e.g., a superconducting tape. Including a nonlinear ferromagnetic material in $\Oa$ does not raise any additional issue. We restrict the analysis to 2D problems with an in-plane magnetic field. With the same procedure as for the \hafOnly, for every iteration $k$ at time step $n$, we obtain the following discrete linear system for the unknowns $\t\in \tspzd(\Gamma_w)$ and $\a\in \aspzd(\Oa)$:
\begin{equation}\label{oneIterationTASystem}
\begin{aligned}
&\volInt{\nu\, \curl \a}{\curl \a'}{\Oa} - \surInt{w\, \curl \t}{\a'}{\Gamma_w} = 0,\\
 -&\surInt{w\, \curl \t'}{\a}{\Gamma_w}-\surInt{\Delta t\, w\, (\partial\e/\partial\j)^{k-1}\ \curl \t}{\curl \t'}{\Gamma_w}\\
& \qquad \qquad \qquad \qquad \qquad\qquad   = \langle \tilde{\vec s}_{\t}, \t' \rangle,
\end{aligned}
\end{equation}
with a right-hand side functional $\tilde{\vec s}_{\t}$ defined by
\begin{align}
&\langle \tilde{\vec s}_{\t}, \t' \rangle = -\surInt{w\,\a_{n-1}}{\curl \t'}{\Gamma_w} + \Delta t \sum_{i\in C} V_i \mathcal{I}_i(\t')\notag \\
&\quad + \surInt{\Delta t\, w \paren{\paren{\rho \mat I - \partial\e/\partial\j} \curl \t}^{k-1}}{\curl \t'}{\Gamma_w}.
\end{align}

In $\aspzd(\Oa)$, we use the same norm as for the \hafOnly,
\begin{align}
\|\a\|_{\aspzd}^2 &= \volInt{\nu_0\ \curl \a}{\curl \a}{\Oa},
\end{align}
and we have $\alpha^\delta \ge \min(\nu/\nu_0)>0$, and $\|a^\delta\| \le \max(\nu/\nu_0)<\infty$, whatever the mesh.

For the discrete inf-sup condition, to avoid the evaluation of a $H^{-1/2}(\Gamma_w)$-norm, we use a mesh-dependent norm, as is common in the discrete setting \cite{bechet2009stable}. We assume a uniform mesh on $\Gamma_w$, for which there exists a $\delta$ and two finite non-zero constants $c_1$ and $c_2$ such that $c_1 \delta \le \delta_e\le c_2 \delta$, $\forall e\in \Gamma_w$, with $\delta_e$ the length of edge $e$. For a given mesh-length $\delta$, we define
\begin{align}\label{eqn_meshdependentnorm}
\|\t\|_{\tspzd}^2 &= \delta \surInt{w\ \Delta t_0\  \rho_0\ \curl \t}{\curl \t}{\Gamma_w},
\end{align}
with $\Delta t_0$ and $\rho_0$ being characteristic time step and resistivity values. The inverse inequality \cite{el2001stability}
\begin{align}
\|\mu\|_{H^{-1/2}(\Gamma)} \ge c\sqrt{\delta} \|\mu \|_{L^2(\Gamma)}, \quad \forall\mu \in H^{-1/2}(\Gamma),
\end{align}
with a finite constant $c$ implies that satisfying the inf-sup test with norm \eqref{eqn_meshdependentnorm} is a necessary condition for stability in terms of norm $\|\cdot \|_{H^{-1/2}(\Gamma_w)}$. In \cite{pitkaranta1979boundary}, the condition is also shown to be sufficient.

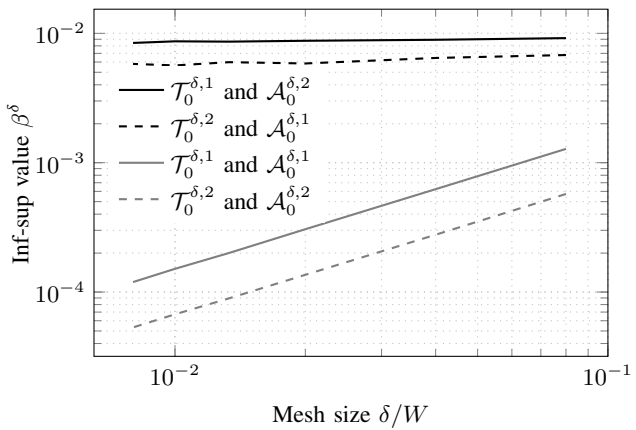
\begin{figure}[h!]
        \centering
	\begin{tikzpicture}[trim axis left, trim axis right][font=\small]
 	\begin{loglogaxis}[
	tick scale binop=\times,
    width=0.95\linewidth,
    height=6.2cm,
    grid = both,
    grid style = dotted,
    xmin=0.0065, 
    xmax=0.1,
	xtick scale label code/.code={},
    xtick={0.01, 0.1, 1},
	ytick scale label code/.code={},
	xlabel={Mesh size $\delta/W$},
    ylabel={Inf-sup value $\beta^\delta$},
    ylabel style={yshift=-0.8em},
    legend style={at={(0.03, 0.61)}, anchor=west, draw=none}
    ]
    \addplot[black, thick] 
    table[x=meshsize,y=beta21]{stabilityMeta/infsup_ta_tapeSuper.txt};
    \addplot[black, dashed, thick] 
    table[x=meshsize,y=beta12]{stabilityMeta/infsup_ta_tapeSuper.txt};
    \addplot[gray, thick] 
    table[x=meshsize,y=beta11]{stabilityMeta/infsup_ta_tapeSuper.txt};
        \addplot[gray, dashed, thick] 
    table[x=meshsize,y=beta22]{stabilityMeta/infsup_ta_tapeSuper.txt};
    \legend{$\tspzdone$ and $\aspzdtwo$, $\tspzdtwo$ and $\aspzdone$, $\tspzdone$ and $\aspzdone$,$\tspzdtwo$ and $\aspzdtwo$}
    \end{loglogaxis}
	\end{tikzpicture}%
	        \caption{Evolution of the inf-sup constant with mesh refinement ($\delta \to 0$) on the simple tape problem ($n=20$, $\jc = 2.5\times 10^{10}$ A/m$^2$). Four cases are considered $\t \in \tspzdi(\Oh)$ and $\a \in \aspzdj(\Oa)$, for $(i,j)\in \{1,2\}\times \{1,2\}$. We observe instabilities when $i= j$. The usual Newton-Raphson scheme with $i=2$, $j=1$ does not converge. Only the case $i=1$, $j=2$ leads to satisfying results.}
        \label{infsup_ta_tapeSuper}
\end{figure}

Fig.~\ref{infsup_ta_tapeSuper} gives the evolution of the inf-sup constant for a sequence of progressively refined meshes, for four choices of function spaces. Analogously to the \hafOnly, it is only when exactly one approximation space is enriched with hierarchical elements that the inf-sup constant is uniformly bounded from below. These choices are good candidates if we want a stable formulation. On the other hand, when choosing $\t \in \tspzdone(\Gamma_w)$ and $\a \in \aspzdone(\Oa)$, or $\t \in \tspzdtwo(\Gamma_w)$ and $\a \in \aspzdtwo(\Oa)$, the test suggests that stability issues may arise. In practice, this is indeed the case, see Fig.~\ref{tape_j}.

We observed that with the choice $\t \in \tspzdtwo(\Gamma_w)$ and $\a \in \aspzdone(\Oa)$, the Newton-Raphson procedure faces convergence troubles. Using a fixed point method does not help either. No satisfying numerical solution has been obtained in the nonlinear case. On the other hand, when considering a linear conductor, no particular issue is encountered and oscillations disappear, as expected from the inf-sup test. This result indicates that this choice for $\t$ and $\a$ is acceptable for linear conductors. We therefore believe that the issue for nonlinear materials is related to the iterative technique rather than to the structure of the saddle-point problem. 

In contrast to that, the choice $\t \in \tspzdone(\Gamma_w)$ and $\a \in \aspzdtwo(\Oa)$ provides good results and no issues have been observed. Our results match the observations in \cite{berrospe2019real}, where the function space for $\a$ is however enriched in the whole $\Oa$ domain, instead of only in the vicinity of $\Gamma_w$.



To conclude, extending to spaces with non-homogeneous essential boundary conditions, we recommend choosing $\t \in \tspdone(\Gamma_w)$ and $\a \in \aspdtwo(\Oa)$. This choice ensures a bounded inf-sup value and does not exhibit any stability issues.


\section{Conclusion}\label{sec_conclusion}

In this work, we presented two coupled finite element formulations. The \haf is efficient for systems containing both superconductors and ferromagnetic materials, whose nonlinear constitutive laws are most efficiently handled by combining different formulations. The so-called \taf is an efficient method for modeling superconducting tapes as surfaces. Two fields are used and coupled on the tapes. We proposed a new derivation of the \taf with global constraints, on either current or voltage for each tape.

Both formulations are mixed on the coupling interfaces and the associated systems of equations take the form of a perturbed saddle point problem. They fit into the classical framework of mixed formulations. It is well known that this problem structure may be exposed to stability issues, e.g., spurious oscillations in the numerical solutions, if function spaces are not chosen correctly. We illustrated the stability issues arising for naive choices of function spaces. We then investigated the formulations stability in the discrete setting, using the classical mixed formulation theory, for several choices of finite element spaces, restricting our study to 2D problems.

The conclusions for both formulations are similar. Using basis functions of different suitable polynomial orders on the coupling interfaces helps to avoid stability issues, whereas with identical orders, the inf-sup value fails to be uniformly bounded above zero. For the \hafOnly, either the space for $\h$, or the space for $\a$ should be enriched, e.g., locally via hierarchical elements on the coupling boundary. For the \tafOnly, the only satisfying configuration consists in using second-order hierarchical elements on the tapes for $\a$ while using first-order elements for $\t$.

Extending to 3D problems would constitute an interesting research topic in further works. Other solutions for stabilizing the problem could also be considered, such as using dual meshes on coupling interfaces for the two fields.

\section*{Appendix}

\subsection{Coupling term in the \haf}

The surface integral to be coupled with the $\a$-field of the \af reads
\begin{align}
\surInt{\e \times \n_{\Oh}}{\h'}{\Gm}.
\end{align}
In this work, $\Gm$ is only placed at the exterior of the conducting domain, or on its boundary. On $\Gm$, the trace $\h'\times \n_{\Oh}$ is therefore locally that of the gradient of a scalar function: $\h'~=~\grad \phi'$ (intersections with possible cut functions are already treated in the global term $V_i \mathcal{I}_i(\h')$). Note that even when $\Gm$ is the boundary of the conducting domain, a scalar potential is introduced on its surface (see Eq.~\eqref{eqn_h_decomposition}). Consequently, we have:
\begin{align}
\surInt{\e \times \n_{\Oh}}{\grad \phi'}{\Gm} =& \surInt{\curl (\phi' \e)}{ \n_{\Oh}}{\Gm} \notag\\
& - \surInt{\phi' \curl\e}{ \n_{\Oh}}{\Gm} .
\end{align}

If $\Gm$ is a closed surface, then the first term in the right-hand side vanishes by Stokes theorem. In the second term, only the curl of $\e$ appears. Because $\e = -\dt \a - \grad v$ in $\Oa$, we have $\curl \e = -\curl (\dt \a)$.

If $\Gm$ is not a closed surface, then $\Gm \cup (\Gamma_e \cap \partial \Oh)\cup (\Gamma_h \cap \partial \Oh)$ is closed. On $(\Gamma_h \cap \partial \Oh)$, $\phi' = 0$, and on $(\Gamma_e \cap \partial \Oh)$, we considered homogeneous natural boundary conditions so $\e\times \n_{\Oh} = \vec 0\Rightarrow \curl \e \cdot \n_{\Oh}=0$. The treatment of non-homogeneous natural boundary conditions is straightforward as well.






\section*{Acknowledgment}

We would like to thank Prof. Barbara Wohlmuth from the Technical University of Munich and Prof. Eric Béchet from the University of Liège for the insightful discussions about mesh-dependent norms. We would also like to thank Prof. Herbert Egger from Technical University of Darmstadt for the fruitful exchanges about saddle-point problems stability analysis.

\bibliographystyle{ieeetr}
\bibliography{paperReferences}

\end{document}